\documentclass[10pt]{article}

\usepackage{dpr_fec,ifthen,dsfont}
\allowdisplaybreaks

\usepackage{isetechreport}
\def\reportyear{17T}
\def\reportno{012}
\def\revisionno{3}
\def\originaldate{February 3, 2018}
\def\revisiondate{\today}
\coraltrue
\cvcrfalse

\begin{document}

\title{A Self-Correcting Variable-Metric Algorithm Framework for Nonsmooth Optimization}

\author{Frank E.~Curtis\thanks{E-mail: frank.e.curtis@gmail.com}}
\affil{Department of Industrial and Systems Engineering, Lehigh University}
\author{Daniel P.~Robinson\thanks{E-mail: daniel.p.robinson@gmail.com}}
\affil{Department of Applied Mathematics and Statistics, Johns Hopkins University}
\author{Baoyu Zhou\thanks{E-mail: baz216@lehigh.edu}}
\affil{Department of Industrial and Systems Engineering, Lehigh University}
\titlepage

\maketitle

%**********
% Abstract
%**********
\begin{abstract}
  An algorithm framework is proposed for minimizing nonsmooth functions.  The framework is variable-metric in that, in each iteration, a step is computed using a symmetric positive definite matrix whose value is updated as in a quasi-Newton scheme.  However, unlike previously proposed variable-metric algorithms for minimizing nonsmooth functions, the framework exploits self-correcting properties made possible through BFGS-type updating.  In so doing, the framework does not overly restrict the manner in which the step computation matrices are updated, yet the scheme is controlled well enough that global convergence guarantees can be established.  The results of numerical experiments for a few algorithms are presented to demonstrate the self-correcting behaviors that are guaranteed by the framework.
\end{abstract}

\bibliographystyle{plain}
% Commands
\newcommand{\inprod}[2]{\langle #1,#2 \rangle}
\newcommand{\matprod}[2]{\langle\!\langle #1,#2 \rangle\!\rangle}
\newcommand{\fb}{f_\diamond}
\newcommand{\gb}{g_\diamond}
\newcommand{\xb}{x_\diamond}
\newcommand{\svano}{\textsc{SVANO}}
\newcommand{\svanobfgs}{\texttt{SVANO-BFGS}}
\newcommand{\svanobfgsstep}{\textsc{SVANO-BFGS-Step}}
\newcommand{\svanobundle}{\texttt{SVANO-Bundle}}
\newcommand{\svanobundlestep}{\textsc{SVANO-Bundle-Step}}
\newcommand{\svanobundlestepnonconvex}{\textsc{SVANO-Bundle-Step-Nonconvex}}
\newcommand{\svanogs}{\texttt{SVANO-GS}}
\newcommand{\bfgs}{\texttt{SVANO-BFGS}}
\newcommand{\bundle}{\texttt{Bundle}}
\newcommand{\gs}{\texttt{GS}}

%*********
% Section
%*********
\section{Introduction}

The purpose of this paper is to present an algorithm framework for solving minimization problems involving nonsmooth (locally Lipschitz) objective functions.  To frame the context and goals of this work, it is worthwhile to recall some history on the design of algorithms over the past few decades.

Practical algorithms for minimizing \emph{smooth} objectives primarily fall between two extremes.  At one extreme are steepest descent methods that only require first-order derivative (i.e., gradient) information.  Such methods have relatively cheap per-iteration costs and can attain a linear rate of convergence to a minimizer.  At the other extreme are Newton methods that require first- and second-order derivative (i.e., respectively, gradient and Hessian) information as well as solving linear systems of equations of dimension equal to the number of variables.  Such methods are relatively expensive, but can attain a quadratic rate of convergence to a minimizer.  For further details, see, e.g., \cite{BazaSherShet06}, \cite{Bert99}, \cite{DennSchn96}, \cite{NoceWrig06}, and \cite{Rusz06}.

One often finds, however, that the most computationally efficient method for a given application does not follow either of these extremes.  That is, one often finds that with only approximate second-order information, and with techniques that avoid expensive computations such as solving linear systems of equations, one can better balance per-iteration costs with per-iteration improvement.

Along these lines, one of the most important developments in the design of smooth optimization algorithms came with the advent of \emph{variable-metric} algorithms in the 1960s; see~\cite{Davi91}.  This class of methods, which includes quasi-Newton methods such as those of the widely successful Broyden-Fletcher-Goldfarb-Shanno (BFGS) variety (see \cite{Broy70}, \cite{Flet70}, \cite{Gold70}, and \cite{Shan70}), often offer an attractive alternative between extremes.  Such methods only require first-order derivative information, can avoid the need to solve linear systems of equations, and yet can offer superlinear convergence rate guarantees; see~\cite{DennMore74}.

When it comes to minimizing \emph{nonsmooth} functions, the array of available algorithms is more varied, and attempts to characterize and compare them run into various challenges.  For example, the ideas underlying steepest descent, quasi-Newton, and Newton methodologies can all be extended for minimizing nonsmooth functions, but practical methods often involve computations beyond obtaining (sub)gradients and solving linear systems---e.g., they include proximal point, cutting plane, gradient sampling, and other methodologies---making the computational trade-offs between methods less straightforward.  In addition, theoretical convergence rates for algorithms become more difficult to prove, meaning that one cannot rely so easily on such characterizations when comparing methods.  For example, despite being introduced decades ago in the 1970s and being one of the most popular classes of methods for convex optimization, convergence rate guarantees for bundle methods have been shown in comparatively fewer articles; for some examples, see~\cite{Robi99}, \cite{Kiwi00}, and \cite{MiffSaga05}.  Moreover, many of these guarantees focus on the rate achieved over the subsequence of accepted (``serious'') steps, not always accounting for the work involved to compute such a step (which might involve some number of intermediate ``null'' steps).  For one exception, see \cite{DuRusz17}.

All of this being said, many have observed that, as in smooth optimization, improved computational trade-offs between per-iteration cost and improvement are often attained by methods that employ both first-order derivative information and approximate second-order information.  In this spirit, this paper proposes a new variable-metric algorithm framework for solving nonsmooth optimization problems.

Variable-metric algorithms for nonsmooth optimization have previously been proposed; see, e.g., \cite{BonnGilbLemaSaga95}, \cite{HiriLema93b}, \cite{Kiwi00}, \cite{Lema82}, \cite{MiffSunQi98}, and \cite{VlceLuks01}.  Broadly speaking, they can be grouped into three categories.  Firstly, there are techniques not built on quasi-Newton-type updating.  An important example in this group is Shor's R-algorithm; see~\cite{Shor85} and the more recent work in \cite{KappKunt00} and \cite{BurkLewiOver08}.  Secondly, there are techniques that attempt to employ quasi-Newton ideas, but only possess convergence guarantees when the updates are restricted to ensure that the resulting Hessian approximations remain sufficiently positive definite and bounded in all iterations; see, e.g., \cite{CurtQue13,CurtQue15}.  Third, there are techniques that employ unadulterated quasi-Newton ideas.  Interestingly, convergence guarantees can be established for unadulterated BFGS in a few specific cases (see~\cite{LewiOver13}), though general guarantees for broad classes of functions remain elusive.

The algorithm framework proposed in this paper falls into the second of the categories in the preceding paragraph, but is unique in that it exploits the \emph{self-correcting properties} of BFGS-type updating.  These properties guarantee that a sufficient number of matrices generated by a BFGS-type updating scheme possess useful properties for ensuring convergence without having to overly restrict the manner in which the updates are performed.  The hope is that the framework proposed in this paper can offer both practical performance gains for various algorithm classes for nonsmooth optimization as well as outline how these useful properties of BFGS-type updating can be incorporated into other algorithms.  The results of our numerical experiments shown in this paper provide evidence that the framework does indeed offer performance gains.  One interesting aspect of our approach, revealed by our numerical experiments, is that it is effective at recognizing when an iterate is nearly stationary.  This is in contrast to the behavior of an unadulterated BFGS approach, which, if employed in practice, might break down (e.g., due to a failed line search) before any stationarity guarantee has been offered.

%************
% Subsection
%************
\subsection{Organization}

In~\S\ref{sec.algorithms}, we state our problem of interest, describe the proposed framework, and discuss at a broad level the types of algorithms that adhere to the framework.  In~\S\ref{sec.convergence}, we discuss the properties of the scaling matrices employed in the framework, then show how these properties can be used to obtain generic convergence guarantees.  In~\S\ref{sec.instances}, we present specific algorithms that adhere to the framework.  The results of numerical experiments are given in \S\ref{sec.numerical}.  Concluding remarks are provided in \S\ref{sec.conclusion}.

%************
% Subsection
%************
\subsection{Notation}

Let $\R{}$ denote the set of real numbers (i.e., scalars), let $\R{}_{\geq0}$ denote the set of nonnegative real numbers, let $\R{}_{>0}$ denote the set of positive real numbers, and let $\N{} := \{1,2,\dots\}$ denote the set of natural numbers.  In addition, for any of these quantities, let a superscript $n \in \N{}$ be used to indicate the $n$-dimensional extension of the set---e.g., let $\R{n}$ denote the set of $n$-dimensional real vectors---and let a superscript $\nbar \times n$ with $(\nbar,n) \in \N{} \times \N{}$ be used to indicate the $\nbar$-by-$n$-dimensional extension of the set---e.g., let $\R{\nbar \times n}$ denote the set of $\nbar$-by-$n$ real matrices.  A vector with all elements equal to 1 is denoted as $\mathds{1}$ and an identity matrix is denoted as $I$, where, in each case, the size of the quantity is determined by the context in which it appears.  With real symmetric matrices $A$ and $B$, let $A$ $\succ$ ($\succeq$) $B$ indicate that $A - B$ is positive definite (semidefinite).  Given a set $\Xcal$, its convex hull is denoted as $\conv\Xcal$.

%*********
% Section
%*********
\section{Problem Statement and Algorithm Framework}\label{sec.algorithms}

In this section, we formally state our optimization problem of interest and our proposed algorithm framework.  We also outline ideas underlying various types of algorithms that adhere to it.

%************
% Subsection
%************
\subsection{Problem Statement}

Our problem of interest is to minimize an objective $f : \R{n} \to \R{}$, i.e., consider the optimization problem
\bequation\tag{P}\label{prob.f}
  \min_{x\in\R{n}}\ f(x).
\eequation
For now, only the following assumption is made about problem~\eqref{prob.f}.

\bassumption\label{ass.f}
  \textit{
  The objective function $f : \R{n} \to \R{}$ in problem~\eqref{prob.f} is bounded below over $\R{n}$, locally Lipschitz on $\R{n}$, and continuously differentiable in an open set $\Dcal$ with full measure in $\R{n}$.
  }
\eassumption

Under this assumption, there exists a scalar $f_{\inf} \in \R{}$ such that
\bequation\label{eq.bounded_below}
  f(x) \geq f_{\inf}\ \ \text{for all}\ \ x \in \R{n},
\eequation
and, for any compact subset $\Bcal$ of $\R{n}$, there exists a constant $L_\Bcal \in \R{}_{>0}$ such that
\bequation\label{eq.Lipschitz}
  |f(x) - f(\xbar)| \leq L_\Bcal \|x - \xbar\|_2\ \ \text{for all}\ \ (x,\xbar) \in \Bcal \times \Bcal.
\eequation
Assumption~\ref{ass.f} does not preclude the possibility that $f$ might have no minimizer, or that it might have many local minimizers.  The goal of our framework is to characterize a family of methods for generating a sequence of iterates that is guaranteed, in the limit, to reveal a stationary point for $f$.  (In other words, our framework is not necessarily intended for global optimization.)  Hence, it is worthwhile to derive stationarity conditions for $f$ that must be satisfied at any local minimizer.

\bremark
  \textit{
  One might be interested in situations when the objective can be unbounded below and/or when it is extended-real-valued.  We discuss such situations in our concluding remarks in~\S\ref{sec.conclusion}.
  }
\eremark

Stationarity conditions for $f$ can be derived following the treatment by~\cite{Clar83}.  (Indeed, many of the following terms are often defined with a ``Clarke'' designation.  We omit this designation for brevity.)  Firstly, the generalized directional derivative of $f$ at $x \in \R{n}$ with respect to $s \in \R{n}$ is given by
\bequationNN
  f^\circ(x;s) = \limsup_{\xbar \to x,\alpha\searrow 0} \frac{f(\xbar + \alpha s) - f(\xbar)}{\alpha}.
\eequationNN
The subdifferential of $f$ at $x$ is then defined as
\bequationNN
  \partial f(x) = \{g \in \R{n} : f^\circ(x;s) \geq g^Ts\ \text{for all}\ s \in \R{n}\}.
\eequationNN
According to Rademacher's theorem, any function $f$ that is locally Lipschitz on~$\R{n}$ is differentiable almost everywhere and its subdifferential at $x$ (see Theorem~2.5.1 in \cite{Clar83}) is given by
\bequationNN
  \partial f(x) = \conv\left\{\lim_{k\to\infty} \nabla f(x_k) : \{x_k\} \to x\ \text{with}\ x_k \in \Dcal\ \text{for all}\ k \in \N{} \right\}.
\eequationNN
For a given $\epsilon \in \R{}_{\geq0}$, the $\epsilon$-subdifferential of $f$ at $x$ (see \cite{Gold77}) is given by
\bequation\label{eq.eps_subdiff}
  \partial_\epsilon f(x) = \conv \partial f(\Bmbb(x,\epsilon)),\ \ \text{where}\ \ \Bmbb(x,\epsilon) := \{\xbar \in \R{n} : \|\xbar - x\|_2 \leq \epsilon\}.
\eequation
A point $x \in \R{n}$ is said to be stationary for $f$ if $0 \in \partial f(x)$ whereas it is merely $\epsilon$-stationary if $0 \in \partial_\epsilon f(x)$.  The following fundamental and widely applicable result, which we attribute to Kiwiel, will be used later.

\blemma\label{lem.kiwiel}
  \textit{
  \emph{(Lemma~3.2(iii), \cite{Kiwi07})}  Let $\{x_k\} \subset \R{n}$ and $\{\epsilon_k\} \subset \R{}_{\geq0}$ be infinite sequences and define $\{\gtilde_k\} \subset \R{n}$ such that $\gtilde_k \in \partial_{\epsilon_k} f(x_k)$ for all $k \in \N{}$.  If, for $x \in \R{n}$, it follows that
  \bequationNN
    \liminf_{k\to\infty}\ \max\{\|x_k - x\|_2,\|\gtilde_k\|_2,\epsilon_k\} = 0,
  \eequationNN
  then $0 \in \partial f(x)$, i.e., the point $x$ is stationary for $f$.
  }
\elemma

%************
% Subsection
%************
\subsection{Algorithm Framework}

The framework that we propose, entitled a \emph{\textbf{S}elf-correcting \textbf{V}ariable-metric \textbf{A}lgorithm for \textbf{N}onsmooth \textbf{O}ptimization}, is stated below as~\ref{alg.svano}.  It consists of two main procedures: $(i)$~Steps~\ref{step.comp_begin}--\ref{step.comp_end}, the computation of a step yielding a reduction in the objective function and $(ii)$~Steps~\ref{step.y}--\ref{step.W}, the computation of quantities used to update a scaling matrix to be used in the step computation procedure in the subsequent iteration.  In~\ref{alg.svano}, these procedures are written in a generic manner so as to allow for flexibility in the choices of various algorithm quantities.  We discuss techniques for performing the step computation procedure next, in \S\ref{sec.step_computations}.  Then, in \S\ref{sec.Hessian_updates}, we motivate the scaling matrix updating strategy.

\floatname{algorithm}{}
\renewcommand\thealgorithm{\svano}
\balgorithm[H]
  \caption{}
  \label{alg.svano}
  \balgorithmic[1]
    \Require A matrix $\Hbar \succ 0$ with smallest (resp.~largest) eigenvalue $\lambda_{\min} \in \R{}_{>0}$ (resp.~$\lambda_{\max} \in \R{}_{>0}$); parameters $\alpha \in (0,1)$, $\eta \in (0,\lambda_{\min}]$, and $\theta \in [\lambda_{\max},\infty)$; a point $x_1 \in \R{n}$; and a positive definite inverse Hessian approximation $W_1 \in \R{n \times n}$.
    \For{\textbf{all} $k \in \N{}$}
      \State Compute, for some $m\in\N{},$ \label{step.comp_begin}
      \bequationNN
        \baligned
        \{x_{k,j}\}_{j=1}^m &\subset \R{n}\ \text{with $x_{k,1} \gets x_k$,} \\
        \{g_{k,j}\}_{j=1}^m &\subset \R{n}\ \text{where $g_{k,j} \in \partial f(x_{k,j})$ for all $(k,j) \in \N{} \times \{1,\dots,m\}$} \\
        \omega_k &\in \R{m}_{\geq0}\ \text{with $\mathds{1}^T\omega_k = 1$,} \\
        \text{and}\ \ \gamma_k &\in \R{n}
        \ealigned
      \eequationNN
      \State such that setting
      \begin{align}
        G_k &\gets \bbmatrix g_{k,1} & \cdots & g_{k,m} \ebmatrix, \label{eq.G} \\
        s_k &\gets -W_k(G_k\omega_k + \gamma_k), \label{eq.s} \\
        \text{and}\ \ x_{k+1} &\gets x_k + s_k \label{eq.x_update}
      \end{align}
      \State yields \label{step.comp_end}
      \bequation\label{eq.sufficient_reduction}
        f(x_{k+1}) \leq f(x_k) - \thalf \alpha (G_k\omega_k + \gamma_k)^TW_k(G_k\omega_k + \gamma_k).
      \eequation
      \State Choose $y_k \in \R{n}$ and compute $\beta_k$ as the smallest value in $[0,1]$ such that \label{step.y}
      \bequation\label{eq.v}
        v_k \gets \beta_k \Hbar s_k + (1 - \beta_k)y_k
      \eequation
      \State yields
      \bequation\label{eq.bounds}
        \eta \leq \frac{s_k^Tv_k}{\|s_k\|_2^2}\ \ \text{and}\ \ \frac{\|v_k\|_2^2}{s_k^Tv_k} \leq \theta,
      \eequation
      \State and then set \label{step.W}
      \bequation\label{eq.W_update}
        W_{k+1} \gets \(I - \frac{v_ks_k^T}{s_k^Tv_k}\)^TW_k\(I - \frac{v_ks_k^T}{s_k^Tv_k}\) + \frac{s_ks_k^T}{s_k^Tv_k}.
      \eequation
    \EndFor
  \ealgorithmic
\ealgorithm

%************
% Subsection
%************
\subsection{Step Computation Techniques}\label{sec.step_computations}

The step computation procedure in~\ref{alg.svano} covers a wide range of nonsmooth optimization methods, including those that employ cutting plane and gradient sampling methodologies using line search and/or trust region techniques.  Given a symmetric positive definite~$W_k$ (i.e., given $W_k \succ 0$), the procedure consists of the selection of a set of points $\{x_{k,j}\}_{j=1}^m$ in the vicinity of (and including) the current iterate $x_k \in \R{n}$ and a set of vectors $\{g_{k,j}\}_{j=1}^m$ where $g_{k,j} \in \partial f(x_{k,j})$ for each pair $(k,j) \in \N{} \times \{1,\dots,m\}$.  (One could instead set $g_{k,j}$ as a convex combination of subgradients of $f$ evaluated at a set of points, as in~\emph{subgradient aggregation}; see, e.g., \cite{Kiwi85}.  However, for simplicity in the algorithm statement, we do not state this option explicitly.)  Following the selection of these vectors, the framework requires a pair $(\omega_k,\gamma_k)$ such that the step~$s_k$ in~\eqref{eq.s} leads to the reduction in $f$ in \eqref{eq.sufficient_reduction}.  The vector $\omega_k$, required to be nonnegative with elements summing to unity, should be viewed as a vector of weights such that, with~$G_k$ defined in \eqref{eq.G}, the step component $G_k\omega_k$ is a convex combination of the elements in~$\{g_{k,j}\}_{j=1}^m$.  The vector~$\gamma_k$ then represents a perturbation of this convex combination, which, e.g., may arise due to the use of line search or trust region methodologies.

Let us make these ideas concrete by describing, for example, how all of these quantities may be derived in a trust region framework.  Suppose that at $x_k \in \R{n}$, a set of points $\{x_{k,j}\}_{j=1}^m$ (with $x_{k,1} \gets x_k$) and vectors $\{g_{k,j}\}_{j=1}^m \subset \R{n}$ are given as in the framework.  In addition, suppose a set of scalars $\{f_{k,j}\}_{j=1}^m$ is given.  (These scalars would typically depend on values of $f$ at $x_k$ and/or $\{x_{k,j}\}_{j=1}^m$; see \S\ref{sec.instances} for more detail.)  Then, a convex piecewise-linear model of $f$ at $x_k$ is given by $l_{k,m} : \R{n} \to \R{}$ defined by
\bequation\label{eq.l}
  l_{k,m}(x) = \max_{j\in\{1,\dots,m\}}\{f_{k,j} + g_{k,j}^T(x - x_{k,j})\}.
\eequation
Also, given $H_k \succ 0$, a convex piecewise-quadratic model of $f$ at $x_k$ is given by $q_{k,m} : \R{n} \to \R{}$ defined by
\bequation\label{eq.q}
  q_{k,m}(x) = l_{k,m}(x) + \thalf (x - x_k)^TH_k(x - x_k).
\eequation
A step toward minimizing $f$ can be defined by the minimizer of $q_{k,m}$ within a region defined by a norm~$\|\cdot\|$ and trust region radius $\delta_k \in \R{}_{>0} \cup \{\infty\}$, i.e., the minimizer of
\bequation\label{prob.primal}
  \min_{x\in\R{n}}\ \ q_{k,m}(x)\ \ \st\ \ x \in \Xcal_k := \{x\in\R{n} : \|x-x_k\| \leq \delta_k\}.
\eequation

Solving \eqref{prob.primal} directly can be challenging due to the nonsmoothness of $l_{k,m}$ (and, hence, of $q_{k,m}$) and due to the presence of the trust region constraint (if $\delta_k < \infty$).  One can reformulate it as the smooth constrained quadratic optimization problem (QP) stated as
\bequation\label{prob.primal_smooth}
  \baligned
    \min_{(x,z)\in\R{n}\times\R{}} &\ z + \thalf (x - x_k)^TH_k(x - x_k) \\
    \st &\ x \in \Xcal_k\ \text{and}\ f_{k,j} + g_{k,j}^T(x - x_{k,j}) \leq z\ \text{for all}\ j \in \{1,\dots,m\},
  \ealigned
\eequation
but even this can be difficult to solve.  Its dual, on the other hand, has properties that might make it easier to solve than~\eqref{prob.primal_smooth}.  Denoting the dual of $\|\cdot\|$ as $\|\cdot\|_*$, the dual of \eqref{prob.primal_smooth} (see Appendix~\ref{app.primal-dual}) is
\bequation\label{prob.dual}
  \sup_{(\omega,\gamma)\in\R{m}_+\times\R{n}} -\thalf(G_k\omega + \gamma)^TW_k(G_k\omega + \gamma) + b_k^T\omega - \delta_k\|\gamma\|_*\ \ \st\ \ \mathds{1}^T\omega = 1,
\eequation
where the vector $b_k\in\R{m}$ has as its $j$th component
\bequation\label{def.b}
  b_{k,j} = f_{k,j} + g_{k,j}^T(x_k - x_{k,j}).
\eequation
The constraints of this dual merely involve an affine equality constraint and lower bounds on some variables.  Therefore, if $\|\cdot\|_*$ is polyhedral---or if $\delta_k = \infty$, in which case the solution $(\omega_k,\gamma_k)$ of \eqref{prob.dual} must have $\gamma_k=0$ (see Appendix~\ref{app.primal-dual})---then it might be more efficient to employ an active-set method to solve~\eqref{prob.dual} than a method for minimizing~\eqref{prob.primal} or one for solving the constrained QP~\eqref{prob.primal_smooth}.  (See \cite{Kiwi86} and \cite{CurtQue13} for more on special-purpose solvers for \eqref{prob.dual} in the case that $(\delta_k,\gamma_k) = (\infty,0)$.)  In any case, if one solves \eqref{prob.dual}, then the solution to~\eqref{prob.primal} can be recovered  as stated as part of the following lemma; for a proof of this result, see Appendix~\ref{app.primal-dual}.

\blemma\label{lem.primal_dual}
  \textit{
  Given the solution $(\omega_k,\gamma_k) \in \R{m}_+ \times \R{n}$ of the dual subproblem~\eqref{prob.dual}, the solution of the primal subproblem~\eqref{prob.primal} is given by $x_k - W_k(G_k\omega_k + \gamma_k)$; hence, $x_{k+1}$ in \eqref{eq.x_update} with $s_k$ given in \eqref{eq.s} is the solution of \eqref{prob.primal}.  In addition, if $f(x_k) \geq l_{k,m}(x_k)$, then one finds that
  \bequation\label{eq.shim}
    f(x_k) - l_{k,m}(x_{k+1}) \geq \thalf (G_k\omega_k + \gamma_k)^TW_k(G_k\omega_k + \gamma_k),
  \eequation
  meaning that, if
  \bequation\label{eq.ared_to_pred}
    f(x_k) - f(x_{k+1}) \geq \alpha (f(x_k) - l_{k,m}(x_{k+1})),
  \eequation
  then \eqref{eq.sufficient_reduction} holds.
  }
\elemma

\bremark
  \textit{
  Observe that the condition $f(x_k) \geq l_{k,m}(x_k)$ in Lemma~\ref{lem.primal_dual} is not restrictive.  If $f$ is convex and $f_{k,j} = f(x_{k,j})$ for all $(k,j) \in \N{} \times \{1,\dots,m\}$, then it is guaranteed since $x \mapsto f_{k,j} + g_{k,j}^T(x - x_{k,j})$ is an affine underestimator of $x \mapsto f(x)$.  On the other hand, if $f$ is nonconvex, then one can ensure the condition using standard techniques.  For example, in a bundle method approach, one can ensure it using standard downshifting ideas; see, e.g.,~\S3 in \cite{SchrZowe92} and further discussion in this paper in~\S\ref{subsec.nonconvex}.  In a gradient sampling method, one chooses $f_{k,j} = f(x_k)$ for all $(k,j) \in \N{} \times \{1,\dots,m\}$, from which one finds that $f(x_k) \geq l_{k,m}(x_k)$ automatically holds.
  }
\eremark

There are also practical benefits of solving the dual~\eqref{prob.dual} when $\{x_{k,j}\}_{j=1}^m$ and corresponding quantities are generated incrementally.  For example, suppose that elements indexed by $j \in \{1,\dots,\mbar\}$ for some $\mbar \in \N{}$ have been generated, but the resulting trial iterate defined as in \eqref{eq.x_update} fails to satisfy \eqref{eq.sufficient_reduction}.  Then, suppose that additional data indexed by $j = \mbar+1$ is generated in some manner to produce the next trial iterate.  The dual subproblem is the same as the previous one, except for the addition of a single dual variable (and corresponding objective and constraint data entries).  The previous optimal dual solution augmented with the new variable initialized to zero represents a feasible solution of the subsequent dual problem, making it an attractive starting point for solving the subsequent dual subproblem.

%************
% Subsection
%************
\subsection{Scaling Matrix Updating Strategy}\label{sec.Hessian_updates}

A critical feature of~\ref{alg.svano} is that each element of the sequence of matrices $\{W_k\}_{k\geq2}$ is set by an update performed during the previous iteration.  The update~\eqref{eq.W_update} has the same form as a standard BFGS update from the smooth optimization literature, and, indeed, the framework is designed to exploit the properties induced by such an update.  However, the framework allows flexibility in the choice of~$y_k$---in theory, any element of $\R{n}$ will suffice---as long as the scalar $\beta_k \in [0,1]$ is chosen such that the bounds in \eqref{eq.bounds} are satisfied; see \cite{Curt16} for the use of this idea for stochastic optimization.  One possible choice for~$y_k$ is the displacement between a subgradient of $f$ at~$x_{k+1}$ with one at~$x_k$, which is natural since this is the choice that can lead to local superlinear convergence guarantees when~$f$ is smooth.  However, given that~\ref{alg.svano} is designed to solve nonsmooth problems, one should not consider this as the only reasonable choice for~$y_k$.  (For example, one might choose the difference between convex combinations of subgradients encountered in the consecutive iterations.)  An important conclusion of our analysis in~\S\ref{sec.convergence} is that the bounds in~\eqref{eq.bounds} are sufficient for ensuring convergence guarantees, and these bounds can be satisfied for any~$y_k$ as long as~$v_k$ is chosen in \eqref{eq.v} with sufficiently large $\beta_k \in [0,1]$.  (The allowed ranges for~$\eta$ and $\theta$ are set so that \eqref{eq.bounds} is well-defined.  Note that $\beta_k = 1$ implies $v_k = \Hbar s_k$, which, due to the ranges for $\eta$ and $\theta$, implies~\eqref{eq.bounds}; hence, $\beta_k \in [0,1]$ always exists such that $v_k$ set by~\eqref{eq.v} satisfies~\eqref{eq.bounds}.)  We discuss a few choices for the matrix $\Hbar \succ 0$ and the sequence~$\{y_k\}$ in \S\ref{sec.numerical}.

As is well known, applying the Sherman-Morrison-Woodbury formula to \eqref{eq.W_update} yields the following updating formula for $\{H_k\}$ where $H_k = W_k^{-1}$ for all $k \in \N{}$:
\bequation\label{eq.H_update}
  H_{k+1} \gets \(I - \frac{s_ks_k^TH_k}{s_k^TH_ks_k}\)^TH_k\(I - \frac{s_ks_k^TH_k}{s_k^TH_ks_k}\) + \frac{v_kv_k^T}{s_k^Tv_k}.
\eequation
Despite \ref{alg.svano} not requiring $\{H_k\}$ explicitly, it is useful to define this sequence.  For one thing, it can be observed from \eqref{eq.s} that the sufficient reduction condition~\eqref{eq.sufficient_reduction} can equivalently be written as
\bequation\label{eq.sufficient_reduction_H}
  f(x_{k+1}) \leq f(x_k) - \thalf \alpha s_k^TH_ks_k,
\eequation
i.e., the condition requires that the reduction in $f$ from $x_k$ to $x_{k+1}$ is proportional to a quadratic function of the step $s_k$, which is a typical requirement for a descent method.  In addition, the properties of the sequence $\{H_k\}$ corresponding to $\{W_k\}$ will be of central importance in~\S\ref{sec.convergence}.

It is worthwhile to mention that \eqref{eq.v} with $\Hbar \equiv H_k$ would reflect a standard \emph{damping} of the BFGS update; see  \cite{Powe78} and \cite{NoceWrig06}.  However, rather than employ an element of the \emph{sequence} $\{H_k\}$ in \eqref{eq.v}, we employ the \emph{fixed} matrix $\Hbar$.  Since this allows us to ensure that \eqref{eq.bounds} holds for the constants $\eta$ and $\theta$ for all $k \in \N{}$, we are able to ensure the self-correcting properties that are central to our convergence analysis.  One cannot maintain such assurances if $H_k$ is used in place of $\Hbar$ in \eqref{eq.v}.  Another alternative would be to employ a sequence $\{\Hbar_k\}$ with eigenvalues uniformly bounded below by~$\lambda_{\min} \in \R{}_{>0}$ and above by $\lambda_{\max} \in \R{}_{>0}$.  That said, for simplicity, let us assume that $\Hbar$ is fixed.

For ease of reference, we refer to $\{H_k\}$ and $\{W_k\}$ as \emph{Hessian approximations} and \emph{inverse Hessian approximations}, respectively.  This terminology should be easy to accept since it is common in the literature on quasi-Newton methods, even for nonsmooth optimization.  However, since~$f$ is nonsmooth, the term ``Hessian'' should be taken loosely as a matrix that approximates changes in the subgradients of~$f$ at nearby points.  See \cite{Clar83} for more information about generalized second derivatives.

%*********
% Section
%*********
\section{Convergence Properties of~\ref{alg.svano}}\label{sec.convergence}

In this section, we explore properties of \emph{any} sequences $\{W_k\}$ and $\{H_k\}$ generated by~\eqref{eq.W_update} and \eqref{eq.H_update}, respectively, then discuss generic convergence properties of the \ref{alg.svano} Framework.  To start, it is immediate that the updates \eqref{eq.W_update} and \eqref{eq.H_update} satisfy \emph{secant-like} equations, namely
\bequation\label{eq.secant}
  W_{k+1}v_k = s_k\ \ \text{and}\ \ H_{k+1}s_k = v_k.
\eequation
One can also derive a geometric interpretation of the updates for the Hessian approximations, revealing that the $k$th update can be viewed as the combination of a \emph{projection} to erase curvature information along~$s_k$---in a sense, temporarily setting $s_k^TH_{k+1}s_k$ to zero---along with a \emph{correction} based on information contained in $v_k$ to yield $s_k^TH_{k+1}s_k = s_k^Tv_k > 0$; see Appendix~\ref{sec.geometric}.  Most importantly for our purposes is that one can show that sequences of such updates result in useful self-correcting properties, which we explore in~\S\ref{subsec:correcting}.  These properties of the Hessian approximations, when cast in terms of the inverse Hessian approximations, yield properties that we use to prove a convergence result in~\S\ref{subsec:svano-convergence}.

%************
% Subsection
%************
\subsection{Self-Correcting Properties of BFGS Updating}\label{subsec:correcting}

It is illustrated in Appendix~\ref{sec.geometric} that the update \eqref{eq.W_update} is a combination of a projection and a correction of the corresponding Hessian approximation.  However, as these updates build upon one another from one iteration to the next, it is important to characterize properties of the resulting matrices and their effects on the computed steps after a \emph{sequence} of updates.  The fact that we show in this subsection is that as long as $v_k$ is chosen to satisfy the two critical inequalities in \eqref{eq.bounds} for all $k \in \N{}$, then despite curvature information along $\linspan(s_k)$ being projected out with the update \eqref{eq.W_update}, the corresponding correction ensures that the sequences of Hessian and inverse Hessian approximations satisfy useful inequalities.

Early work on the convergence of quasi-Newton methods by \cite{Powe76} and others (see, e.g.,  \cite{ByrdNoceYuan87}, \cite{Ritt79}, \cite{Ritt81}, and \cite{Wern78}) involved analyses that bound the growth of the traces and the determinants of $\{H_k\}$.  In what follows, we follow the work in \cite{ByrdNoce89} involving a streamlined approach in which one bounds the growth of a function defined by a combination of these quantities; see also the summary provided in \cite{NoceWrig06}.

Given $H \succ 0$, consider $\psi : \R{n\times n} \to \R{}$ defined by $\psi(H) = \trace(H) - \ln(\det(H))$.  It can be shown that $\psi(H)$ is positive (in fact, at least $n$) and represents a measure of closeness between $H$ and the identity matrix $I$ (for which $\psi(I) = n$); in particular,~$\psi(H)$ is an upper bound for the natural logarithm of the condition number of $H$. In addition, the update \eqref{eq.H_update} implies that, for all $k \in \N{}$, one has
\bsubequations\label{eq.tr_det_updates}
  \begin{align}
    \trace(H_{k+1}) &= \trace(H_k) - \frac{\|H_ks_k\|_2^2}{s_k^TH_ks_k} + \frac{\|v_k\|_2^2}{s_k^Tv_k} \label{eq.tr_update} \\ \text{and (see \cite{Pear69})}\ \ 
    \det(H_{k+1}) &= \det(H_k) \(\frac{s_k^Tv_k}{s_k^TH_ks_k}\), \label{eq.det_update}
  \end{align}
\esubequations
with which one can explicitly relate $\psi(H_{k+1})$ and $\psi(H_k)$.  Specifically, assuming that $H_k \succ 0$ and the iterate displacement satisfies $s_k \neq 0$, then by defining
\bequation\label{eq.cos}
  \cos\phi_k := \frac{s_k^TH_ks_k}{\|s_k\|_2\|H_ks_k\|_2}\ \ \text{and}\ \ \iota_k := \frac{s_k^TH_ks_k}{\|s_k\|_2^2}
\eequation
it follows from \eqref{eq.tr_det_updates} that
\bequation\label{eq.gamma}
  \psi(H_{k+1}) = \psi(H_k) + \underbrace{\frac{\|v_k\|_2^2}{s_k^Tv_k} - 1 - \ln\(\frac{s_k^Tv_k}{\|s_k\|_2^2}\)}_{\in\R{}} + \underbrace{\ln(\cos^2\phi_k)}_{\leq 0} + \underbrace{\(1 - \frac{\iota_k}{\cos^2\phi_k} + \ln\(\frac{\iota_k}{\cos^2\phi_k}\)\)}_{\leq 0}. 
\eequation
Nonpositivity of the latter two terms is easily verified; see Appendix~\ref{app.self-correcting}.

By restricting the growth of $\psi$ over $\{H_k\}$ and noting that there must exist certain iterations in which the latter terms in \eqref{eq.gamma} are not too negative, one can prove the following theorem showing \emph{self-correcting properties} of the update \eqref{eq.H_update}.  For completeness, we provide a proof of this theorem in Appendix~\ref{app.self-correcting}; see also Theorem~2.1 in \cite{ByrdNoce89}.

\btheorem\label{th.self-correction}
  \textit{
  Let $\{H_k\}$ satisfy~\eqref{eq.H_update} and suppose that there exist $(\eta,\theta) \in \R{}_{>0} \times \R{}_{>0}$ such that \eqref{eq.bounds} holds for all $k \in \N{}$.  Then, for any $p \in (0,1)$, there exist constants $(\kappa,\sigma,\mu) \in \R{}_{>0} \times \R{}_{>0} \times \R{}_{>0}$ such that, for any $K \in \{2,3,\dots\}$, the following hold for at least $\lceil pK \rceil$ values of $k \in \{1,\dots,K\}$:
  \bequation\label{eq.good_bounds_H}
    \kappa \leq \frac{s_k^TH_ks_k}{\|s_k\|_2\|H_ks_k\|_2}\ \ \text{and}\ \ \sigma \leq \frac{\|H_ks_k\|_2}{\|s_k\|_2} \leq \mu.
  \eequation
  }
\etheorem

This theorem leads to the following corollary about the inverse approximations.

\bcorollary\label{cor.self-correction}
  \textit{
  Let $\{W_k\}$ satisfy~\eqref{eq.W_update} and suppose that there exist $(\eta,\theta) \in \R{}_{>0} \times \R{}_{>0}$ such that~\eqref{eq.bounds} holds for all $k \in \N{}$.  Then, for any $p \in (0,1)$, there exist constants $(\nu,\xi) \in \R{}_{>0} \times \R{}_{>0}$ such that, for any $K \in \{2,3,\dots\}$, the following hold for at least $\lceil pK \rceil$ values of $k \in \{1,\dots,K\}$:
  \bequation\label{eq.good_bounds}
    \begin{aligned}  
      \nu \|G_k\omega_k+\gamma_k\|_2^2 
      &\leq (G_k\omega_k+\gamma_k)^TW_k(G_k\omega_k+\gamma_k)\ \ \text{and} \\
      \|W_k(G_k\omega_k+\gamma_k)\|_2^2
      &\leq \xi\|G_k\omega_k+\gamma_k\|_2^2.
    \end{aligned}
  \eequation
  }
\ecorollary
\bproof
  Since the elements of $\{W_k\}$ satisfy~\eqref{eq.W_update}, it follows that the elements of $\{H_k\} = \{W_k^{-1}\}$ satisfy~\eqref{eq.H_update}.  Hence, the conditions of Theorem~\ref{th.self-correction} hold, meaning that the conclusions of Theorem~\ref{th.self-correction} hold.  Then, with \eqref{eq.s}, the inequalities in \eqref{eq.good_bounds_H} can be rewritten using the notation $\gbar_k := G_k \omega_k + \gamma_k$ as
  \bequation\label{eq:gbar-bounds}
    \kappa \leq \frac{\gbar_k^TW_k\gbar_k}{\|W_k\gbar_k\|_2\|\gbar_k\|_2}\ \ \text{and}\ \ \sigma \leq \frac{\|\gbar_k\|_2}{\|W_k\gbar_k\|_2} \leq \mu.
  \eequation
  From the first and third of the inequalities in~\eqref{eq:gbar-bounds}, it follows that
  \bequationNN
    \gbar_k^TW_k\gbar_k \geq \kappa\|W_k\gbar_k\|_2\|\gbar_k\|_2 \geq (\kappa/\mu) \|\gbar_k\|_2^2,
  \eequationNN
  so that the first inequality in~\eqref{eq.good_bounds} holds with  $\nu := \kappa/\mu$.  Meanwhile, from the second inequality in~\eqref{eq:gbar-bounds}, it follows that $\|W_k\gbar_k\|_2^2 \leq \sigma^{-2}\|\gbar_k\|_2^2$, so that the second inequality in~\eqref{eq.good_bounds} holds with $\xi := \sigma^{-2}$.
\eproof

The role played by $p$ in Theorem~\ref{th.self-correction} and Corollary~\ref{cor.self-correction} can be understood as follows.  For any given $p \in (0,1)$, the results show that at least a fraction $p$ of iterations---i.e., at least $\lceil pK \rceil$ out of any $K$---will involve \emph{good} approximations in the sense that there exist constants such that \eqref{eq.good_bounds_H} and \eqref{eq.good_bounds} hold.  Since one can consider~$p$ to be arbitrarily close to 1, one can claim that nearly all iterations involve good approximations.  That said, the constants might be worse for $p$ closer to 1; e.g., the closer $p$ is to 1, the smaller might be $\nu \in \R{}_{>0}$ and the larger might be $\xi \in \R{}_{>0}$.  These constants also depend on the values $(\eta,\theta)$ employed in the algorithm; we remark on the empirical influence of these values in \S\ref{sec.numerical}.

%************
% Subsection
%************
\subsection{Convergence Results for~\ref{alg.svano}}\label{subsec:svano-convergence}

In this subsection, we provide a couple fundamental convergence results for~\ref{alg.svano}.  Here, our goal is to prove generic results that are useful in various circumstances.  In~\S\ref{sec.instances}, we present a few algorithm instances that fall under the~\ref{alg.svano} Framework.  For concreteness in that section, we prove a complete convergence theory for one of those instances using the results provided here.

Our first result represents a fundamental component of the convergence theory for any algorithm that falls under the \ref{alg.svano} Framework.  In particular, it shows that there exists an infinite subsequence of iterations in which the required decreases in the objective~$f$ guarantee that subsequences of $\{G_k\omega_k + \gamma_k\}$ and $\{s_k\}$ vanish.  The proof reveals that the self-correcting properties of the inverse Hessian approximation scheme are critical.  In particular, since at least a fraction of the approximations are \emph{good}, the reductions in $f$ force the right-hand sides in \eqref{eq.good_bounds} to vanish over a subsequence of iterations, which in turn force the left-hand sides in \eqref{eq.good_bounds} to vanish over the same subsequence of iterations.

\btheorem\label{th.svano}
  \textit{
  There exists an infinite index set $\Kcal \subseteq \N{}$ such that the sequences $\{(G_k,\omega_k,\gamma_k)\}$ and~$\{s_k\}$ computed by~\ref{alg.svano} respectively satisfy
  \bequation\label{eq.svano_limits}
    \lim_{k\in\Kcal,k\to\infty} \|G_k\omega_k+\gamma_k\|_2 = 0\ \ \text{and}\ \ 
    \lim_{k\in\Kcal,k\to\infty} \|s_k\|_2 = 0.
  \eequation
  For example, for a given $p \in (0,1)$, this set $\Kcal$ may at least include the infinite set of indices for which Corollary~\ref{cor.self-correction} guarantees the existence of $(\nu,\xi) \in \R{}_{>0} \times \R{}_{>0}$ such that \eqref{eq.good_bounds} holds for all $k \in \Kcal$.
  }
\etheorem
\bproof
  It follows by \eqref{eq.sufficient_reduction} that, for all $k \in \N{}$, one has
  \bequation\label{eq.f_reduction}
    f(x_{k+1}) \leq f(x_k) - \thalf \alpha (G_k\omega_k+\gamma_k)^TW_k(G_k\omega_k+\gamma_k).
  \eequation
  For a given $p \in (0,1)$, let $\Kcal \subseteq \N{}$ be the infinite set of indices for which Corollary~\ref{cor.self-correction} guarantees the existence of $(\nu,\xi) \in \R{}_{>0} \times \R{}_{>0}$ such that \eqref{eq.good_bounds} holds for all $k \in \Kcal$.  Then, for all $k \in \Kcal$, it follows from~\eqref{eq.f_reduction} and the first inequality in \eqref{eq.good_bounds} that
  \bequationNN
    f(x_{k+1}) \leq f(x_k) - \thalf \nu \alpha \|G_k\omega_k+\gamma_k\|_2^2.
  \eequationNN
  Since $f$ is bounded below (see~\eqref{eq.bounded_below}) and monotonically decreasing, the first limit in~\eqref{eq.svano_limits} holds.  Combining it with the second inequality in \eqref{eq.good_bounds} and $s_k$ from~\eqref{eq.s}, the second limit in \eqref{eq.svano_limits} follows.
\eproof

The conclusions of Theorem~\ref{th.svano} are not entirely consequential in their own right.  However, the theorem is fundamental in that it can be used to show that if the columns of $G_k$ correspond to (convex combinations of) subgradients of $f$ evaluated at points in the vicinity of $x_k$ for all $k \in \N{}$, then, as long as the vanishing of $\{G_k\omega_k + \gamma_k\}$ implies the vanishing of $\{G_k\omega_k\}$ (at least over a subsequence), the first limit in \eqref{eq.svano_limits} must mean that a stationary point of $f$ is revealed by a subsequence of the iterates.  In order to have a formal result along these lines to which we can refer later, we state the following theorem, which may be viewed as a more practical version of Lemma~\ref{lem.kiwiel}.

\btheorem\label{th.stationary}
  \textit{
  Suppose that there exists an infinite index set $\Kcal' \subseteq \N{}$ such that
  \bequation\label{eq.Gw_lim}
    \lim_{k\in\Kcal',k\to\infty} \|G_k\omega_k\|_2 = 0.
  \eequation
  In addition, suppose that for all $k \in \Kcal'$ there exists $\epsilon_k \in \R{}_{\geq0}$ such that, for all $j \in \{1,\dots,m\}$, the vector $g_{k,j}$ is a subgradient (or convex combination of subgradients) of $f$ evaluated at a finite subset of $\Bmbb(x_k,\epsilon_k)$ as defined in~\eqref{eq.eps_subdiff}.  Then, if for some $x \in \R{n}$ one has
  \bequation\label{eq.K_liminf}
    \liminf_{k\in\Kcal',k\to\infty} \max\{\|x_k - x\|_2,\epsilon_k\} = 0,
  \eequation
  then $0 \in \partial f(x)$, i.e., the limit point $x$ of $\{x_k\}_{k\in\Kcal'}$ is stationary for $f$.
  }
\etheorem
\bproof
Under the stated conditions, it follows that $g_{k,j}\in\partial_{\epsilon_k} f(x_k)$ for all $(k,j) \in \Kcal' \times \{1,\dots m\}$. Since $\partial_{\epsilon_k} f(x_k)$ is convex by definition and $\gtilde_k := G_k\omega_k$ is a convex combination of $\{g_{k,j}\}_{j=1}^m$, it follows that $\gtilde_k\in\partial_{\epsilon_k} f(x_k)$ for all $k\in\Kcal'$.  Hence, with \eqref{eq.Gw_lim} and \eqref{eq.K_liminf}, the result follows from  Lemma~\ref{lem.kiwiel}.
\eproof

This theorem reveals useful consequences of the first limit in \eqref{eq.svano_limits}, insofar as this limit might be useful toward proving \eqref{eq.Gw_lim}.  But what, if any, are useful consequences of the second limit in \eqref{eq.svano_limits}?  This depends on the instance of \ref{alg.svano} of interest.  For instances where, perhaps under additional assumptions about the objective function $f$, one can guarantee the existence of $\{\epsilon_k\}_{k\in\Kcal'}$ and $x \in \R{n}$ such that \eqref{eq.K_liminf} holds, the second limit is not of great interest in itself.  However, for other instances, having the subsequence of step norms converging to zero helps to ensure that \eqref{eq.K_liminf} holds for some $\{\epsilon_k\}_{k\in\Kcal'}$ and $x \in \R{n}$.  This is demonstrated for our specific instance in \S\ref{subsec.bundle}.

%*********
% Section
%*********
\section{Instances of~\ref{alg.svano}}\label{sec.instances}

Our goal in this section is to present instances of~\ref{alg.svano} that yield the convergence guarantee in Theorem~\ref{th.stationary}. Specifically, we consider choices for computing the quantities in Step~\ref{step.comp_begin} of~\ref{alg.svano}, and for these quantities establish that~\eqref{eq.sufficient_reduction}, \eqref{eq.Gw_lim}, and~\eqref{eq.K_liminf} hold.   However, before doing so, it is instructive to show how, if one tries to fit a classical BFGS strategy into the framework, certain behaviors might cause the method to falter.  This helps to motivate the more involved strategies that we present.

We stress that we are not attempting to claim that one cannot prove convergence for an algorithm that makes use of a classical BFGS strategy.  We merely hope to illustrate the gaps in the analysis that might arise if one were to use such a strategy and try to employ Theorem~\ref{th.stationary} to prove convergence.

%************
% Subsection
%************
\subsection{Classical BFGS Method}\label{subsec.bfgs}

Let us follow \cite{LewiOver13} and discuss a BFGS method with a weak Wolfe line search.  A description of a step computation procedure for such an algorithm (written, in contrast to \cite{LewiOver13}, with subgradients instead of gradients) is presented as \ref{alg.bfgs}.  This procedure should be viewed as an instance of the step computation written as Steps~\ref{step.comp_begin}--\ref{step.comp_end} in \ref{alg.svano}.  Rather than delineate the details of a weak Wolfe line search, we direct the reader to Algorithm~4.6 in \cite{LewiOver13} and note that our required condition \eqref{eq.suff_dec} corresponds to ``$c_1$''$ = \thalf\alpha$, ``$s$''$ = g_{k,1}^T\stilde_k$, and ``$t$''$ = \tilde\alpha_k^2$, which with the subsequent choices for $s_k$ and $x_{k+1}$ implies that~\eqref{eq.suff_dec} yields \eqref{eq.sufficient_reduction} with $G_k\omega_k = g_{k,1}$ and $\gamma_k = (\tilde\alpha_k - 1)g_{k,1}$.  When the line search fails to produce a sufficiently large stepsize, the algorithm \emph{breaks down}; see \S6.1 in \cite{LewiOver13}.

\floatname{algorithm}{}
\renewcommand\thealgorithm{\svanobfgsstep}
\balgorithm[t]
  \caption{}
  \label{alg.bfgs}
  \balgorithmic[1]
    \Require A minimum stepsize parameter $\tilde\alpha_{\min} \in \R{}_{>0}$.
    \State Set $x_{k,1} \gets x_k$ and $g_{k,1} \in \partial f(x_k)$.
    \State Set $G_k \gets \bbmatrix g_{k,1} \ebmatrix$, $\omega_k \gets 1$, and $\stilde_k \gets -W_kG_k\omega_k$.
    \State Run a weak Wolfe line search \cite[Alg.~4.6]{LewiOver13} from $x_k$ along $\stilde_k$ to set
    \bequation\label{eq.suff_dec}
      \tilde\alpha_k \geq \tilde\alpha_{\min}\ \ \text{with}\ \ f(x_k + \tilde\alpha_k \stilde_k) \leq f(x_k) + \thalf \alpha \tilde\alpha_k^2 g_{k,1}^T\stilde_k,
    \eequation
    or \textbf{terminate} (i.e., \emph{break down}) if no such $\tilde\alpha_k$ is found within an iteration limit.
    \State Set $s_k \gets \tilde\alpha_k\stilde_k$ (meaning $s_k \gets -W_k(G_k\omega_k + \gamma_k)$ with $\gamma_k = (\tilde\alpha_k - 1)g_{k,1}$).
    \State Set $x_{k+1} \gets x_k + s_k$.
    \State \textbf{return} $(s_k,x_{k+1})$ to Step~\ref{step.y} in~\ref{alg.svano}.
  \ealgorithmic
\ealgorithm

The issues that may arise for this algorithm all relate to the line search.  If for some $k \in \N{}$ the function~$f$ is not differentiable at $x_k$, then $\stilde_k \gets -W_kG_k\omega_k = -W_kg_{k,1}$ might not be a descent direction for $f$ from $x_k$.  For this and other reasons (see \S4 in \cite{LewiOver13}), the line search might not be able to produce a stepsize within a prescribed iteration limit such that~\eqref{eq.suff_dec} holds.  In such cases, one cannot guarantee the conditions of Theorem~\ref{th.stationary}, namely, \eqref{eq.Gw_lim}, since vanishing of the sequence $\{G_k\omega_k+\gamma_k\} = \{g_{k,1} + (\tilde\alpha_k - 1)g_{k,1}\}$ might not correspond to vanishing of (a subsequence of) $\{G_k\omega_k\} = \{g_{k,1}\}$ due to the large perturbations $\{\gamma_k\} = \{(\tilde\alpha_k - 1)g_{k,1}\}$.  One can imagine various heuristics such that, if the line search would otherwise break down, one might replace $\stilde_k$ with $-W_kG_k\omega_k$ for some matrix of subgradients $G_k$ evaluated at points in $\Bmbb(x_k,\epsilon_k)$ for some $\epsilon_k \in \R{}_{>0}$ and some nonnegative weight vector~$\omega_k$ that sums to unity.  However, it is a nontrivial task to determine such quantities to ensure that the weak Wolfe line search will be guaranteed to return a stepsize above a prescribed positive threshold.

All of this being said, \emph{if} \ref{alg.bfgs} yields $\tilde\alpha_k \geq \tilde\alpha_{\min}$ for all $k \in \N{}$ (which occurs often in practice for sufficiently small, but reasonable values of~$\tilde\alpha_{\min}$, at least until the algorithm has nearly approached a stationary point), then the resulting instance of \ref{alg.svano} attains some of the guarantees in~\S\ref{subsec:svano-convergence}.  In particular, with $\tilde\alpha_k \geq \tilde\alpha_{\min}$ for all $k \in \N{}$, the conditions of Theorem~\ref{th.stationary} hold, and if $g_{k,1} \in \partial f(x_k)$ or some heuristic is used to set $\stilde_k \gets -W_kG_k\omega_k$ as described in the previous paragraph, then the supposition about $\{g_{k,j}\}_{j=1}^m$ in Theorem~\ref{th.stationary} also holds.  Consequently, \emph{if} a subsequence of iterates converges to a limit and the subgradients employed in the step computation are evaluated at points in narrowing neighborhoods of each iterate, then a stationary point will be revealed by $\{x_k\}$.

Other classical BFGS variants can be derived that employ a trust region mechanism instead of a line search.  However, the issues for such methods would be similar to those described above: \emph{if} one finds that a successful step is taken sufficiently often when the trust region radius is above a positive threshold (say, proportional to $\|G_k\omega_k\|_2$), then the algorithm has the potential to converge.  However, it is nontrivial to design an algorithm that maintains the spirit of a basic quasi-Newton algorithm and ensures that this occurs under loose assumptions on $f$.

%************
% Subsection
%************
\subsection{A Bundle Trust Region Method for Convex Minimization}\label{subsec.bundle}

Bundle methods are an extremely popular class of algorithms for solving nonsmooth optimization problems.  Modern variants of bundle methods have guarantees for solving both convex (see, e.g., \cite{HiriLema93b}, \cite{Kiwi85b}, \cite{MiffSaga05}, \cite{LemaNemiNest95}, and \cite{Rusz06}) and nonconvex (see, e.g.,  \cite{ApkaNollProt08}, \cite{HaarMietMaek07}, \cite{HareSaga10}, \cite{Kiwi85}, \cite{Kiwi96}, \cite{LuksVlce98}, \cite{Miff82}, and \cite{SchrZowe92}) problems.  Those for solving convex problems are based on combining cutting plane and proximal point methodologies, whereas those for solving nonconvex problems often employ the same ideas with modifications involving ``downshifting'' and ``tilting'' of the cutting planes.  We refer the reader to the references above for further information and details.

In this subsection, we present an instance of \ref{alg.svano} using bundle method ideas for minimizing convex $f$.  Our instance is described by specifying an algorithm for the step computation procedure written generically in Steps~\ref{step.comp_begin}--\ref{step.comp_end} in \ref{alg.svano}.  In particular, see~\ref{alg.bundle} below.  Through an inner loop, the procedure computes trial steps through successive subproblem solves until one is computed satisfying a descent condition, shown in our analysis below to imply that~\eqref{eq.sufficient_reduction} holds.  (This inner loop causes ``null'' steps to occur until a ``serious'' step is computed, which causes the inner loop to terminate.)  Each iteration of the inner loop involves solving a subproblem of the form \eqref{prob.primal} by solving the dual~\eqref{prob.dual} for $(\omega_{k,m},\gamma_{k,m}) \in \R{m}_+ \times \R{n}$ and recovering the primal solution as described in~\S\ref{sec.step_computations}.  (We now include a second subscript on the solution vectors since they change as the algorithm iterates over $m \in \N{}$.)  The ``bundles'' employed in the loop are the tuples $\{(x_{k,j},f_{k,j},g_{k,j})\}_{j=1}^m$ where $\{f_{k,j}\}_{j=1}^m$ and $\{g_{k,j}\}_{j=1}^m$ are objective values and subgradients evaluated at $\{x_{k,j}\}_{j=1}^m$.  As in a standard bundle method, these include elements of bundles computed in previous ``outer'' iterations, but, for simplicity and since it is not required for convergence, we do not state this in the algorithm explicitly.

\floatname{algorithm}{}
\renewcommand\thealgorithm{\svanobundlestep}
\balgorithm[ht]
  \caption{}
  \label{alg.bundle}
  \balgorithmic[1]
    \State Set $x_{k,1} \gets x_k$, $f_{k,1} \gets f(x_{k,1})$, and $g_{k,1} \in \partial f(x_{k,1})$.
    \For{\textbf{all} $m \in \N{}$}
      \State Set $G_{k,m} \gets \bbmatrix g_{k,1} & \cdots & g_{k,m} \ebmatrix$, then compute $(\omega_{k,m},\gamma_{k,m})$ by solving~\eqref{prob.dual}.
      \State \label{step.bundle.dual} Set $x_{k,m+1} \gets x_k - W_k(G_{k,m}\omega_{k,m} + \gamma_{k,m})$.
      \State Set $f_{k,m+1} \gets f(x_{k,m+1})$ and $g_{k,m+1} \in \partial f(x_{k,m+1})$.
      \State \textbf{if} $l_{k,m}(x_{k,m+1}) = f(x_k)$ \textbf{then} \textbf{terminate} since $x_k$ is stationary for $f$.\label{step.termination_check}
      \If{$f(x_k) - f_{k,m+1} \geq \alpha (f(x_k) - l_{k,m}(x_{k,m+1}))$} \label{step.decrease}
         \State Set $s_k \gets -W_k(G_{k,m}\omega_{k,m} + \gamma_{k,m})$ and $x_{k+1} \gets x_{k,m+1}$.
         \State \Return $(s_k, x_{k+1})$ to Step~\ref{step.y} in~\ref{alg.svano}.
      \EndIf
    \EndFor
  \ealgorithmic
\ealgorithm

For concreteness in order to demonstrate a complete analysis for an algorithm that falls under the \ref{alg.svano} Framework, we show that a standard convergence analysis for a bundle method can be adapted in order to prove a convergence guarantee for \ref{alg.svano} with step computations using \ref{alg.bundle}.  We present an analysis that proceeds in two stages.  Firstly, it is argued that, for any $k \in \N{}$ in which the iterate $x_k$ is suboptimal, the inner loop will terminate finitely, ensuring that the algorithm is well-defined in the sense that it will either reach a minimizer of $f$ in a finite number of iterations or generate an infinite sequence of outer iterates. For this, we borrow from results that are common in the literature on bundle methods. In particular, our first few lemmas---for which complete proofs can be found in Appendix~\ref{app.bundle}---follow the treatment in Chapter~7.4 of \cite{Rusz06}, which in turn borrows from \cite{BonnGilbLemaSaga06} and \cite{HiriLema93b}. (Motivation and analyses for algorithms that combine bundle and trust-region-like ideas go back further as well; see, e.g., \S4 of \cite{Kiwi90}, the article by \cite{SchrZowe92} and references therein, and the book \cite{MakeNeit92}.)  Secondly, it is argued that each accepted step yields a sufficient reduction in~$f$ such that any limit point of the outer iteration sequence is a solution of problem \eqref{prob.f}.  For this, it is important to recognize that standard results cannot readily be applied since the inverse Hessian updating scheme in \ref{alg.svano} does not guarantee that uniformly positive definite and bounded approximations will be generated in all iterations. That said, we are still able to establish a meaningful result due to the critical self-correcting properties of the updating scheme established in~\S\ref{subsec:correcting}.

Toward proving that~\ref{alg.bundle} is well-defined, one may use a type of \emph{Moreau-Yosida regularization} function of~$f$ corresponding, for a given $k \in \N{}$, to the symmetric positive-definite~$H_k$ and trust region $\Xcal_k$ (recall~\eqref{prob.primal}); specifically, consider the function $f_{H_k,\Xcal_k} : \R{n} \to \R{}$ defined by
\bequation\label{eq.MY}
  f_{H_k,\Xcal_k}(\xbar) = \min_{x\in\Xcal_k} f(x) + \thalf(x-\xbar)^TH_k(x-\xbar).
\eequation
This function provides a mechanism for quantifying the separation between $f$ and the models~$l_{k,m}$ and~$q_{k,m}$ defined in \eqref{eq.l} and \eqref{eq.q}, respectively.  To start, the following lemma states that if $x_k$ is not a minimizer of $f$, then this Moreau-Yosida regularization function's value at $x_k$ is strictly less than the objective function value at~$x_k$.  The proof of the result is based simply on the existence of another point in~$\R{n}$ that yields a better objective function value than does~$x_k$.  (Again, the proof of this result and those of Lemma~\ref{lem.MY_bound_q} and Lemma~\ref{lem.bundle.descent} can be found in Appendix~\ref{app.bundle}.)

\blemma\label{lem.MY_bound_f}
  \textit{
  Suppose $f$ is convex. For any $k \in \N{}$, if $x_k$ does not minimize $f$, then $f_{H_k,\Xcal_k}(x_k) < f(x_k)$.
  }
\elemma

The next lemma shows that the Moreau-Yosida regularization function offers an upper bound for the piecewise-linear and piecewise-quadratic model values corresponding to the optimal solution of \eqref{prob.primal}.

\blemma\label{lem.MY_bound_q}
  \textit{
  Suppose $f$ is convex.  For any $(k,m) \in \N{} \times \N{}$, the value of~$l_{k,m}$ evaluated at $x_{k,m+1}$ is bounded above by the optimal value of~\eqref{prob.primal}, which, in turn, is bounded above by the Moreau-Yosida regularization function $f_{H_k,\Xcal_k}$ evaluated at $x_k$; i.e.,
  \bequation\label{eq.lqf}
    l_{k,m}(x_{k,m+1}) \leq q_{k,m}(x_{k,m+1}) \leq f_{H_k,\Xcal_k}(x_k).
  \eequation
  }
\elemma

The following lemma now shows that the inner loop is well-defined in that if a minimizer of~$f$ has not yet been obtained, then a point satisfying the condition in Step~\ref{step.decrease} will be computed.

\blemma\label{lem.bundle.descent}
  \textit{
  Suppose $f$ is convex.  For any $k \in \N{}$, if $x_k$ is not a minimizer of $f$ and $\delta_k \in \R{}_{>0}$ $($recall \eqref{prob.primal}$)$, then \ref{alg.bundle} terminates.
  }
\elemma

Overall, unless the algorithm lands on a point that is optimal for~$f$, the step computation procedure terminates for any $k \in \N{}$.  Since the condition in Step~\ref{step.decrease} implies \eqref{eq.ared_to_pred}, it follows from Lemma~\ref{lem.primal_dual} that the computed step satisfies~\eqref{eq.sufficient_reduction}, meaning that the overall algorithm is well posed.  Thus, all that remains is to show that this instance of \ref{alg.svano} satisfies the remaining assumptions of Theorem~\ref{th.stationary}.  This is done in the following theorem when one introduces a particular strategy for updating the trust region.

\btheorem\label{th.bundle}
  \textit{
  Suppose $f$ is convex.  Consider the \ref{alg.svano} framework in which
  \bitemize
    \item step computations are performed using \ref{alg.bundle}, and
    \item with $\delta_1 \in (0,\infty)$, $\tau \in (0,1)$, and $(\upsilon_1,\upsilon_2,\upsilon_3) \in \R{}_{>0} \times \R{}_{>0} \times \R{}_{>0}$, one sets at the end of each iteration $k \in \N{}$ the next trust region radius as
  \eitemize
  \bequation\label{eq.delta_update}
    \delta_{k+1} \gets \bcases \tau \delta_k & \text{if $\max\{\upsilon_1\|G_k\omega_k + \gamma_k\|_2,\upsilon_2\|s_k\|_2,\upsilon_3\|G_k\omega_k\|_2\} \leq \delta_k$} \\ \delta_k & \text{otherwise.} \ecases
  \eequation
  Then, either the algorithm lands on a minimizer of $f$ in a finite number of iterations or $\{\delta_k\} \searrow 0$ and, with $\Kcal' \subseteq \N{}$ defined as the infinite index set such that $\delta_{k+1} \gets \tau\delta_k$ for all $k \in \Kcal'$, one finds
  \bequation\label{eq.almost_there}
    \lim_{k\in\Kcal',k\to\infty} \|G_k\omega_k\|_2 = 0
  \eequation
  with any limit point of $\{x_k\}_{k\in\Kcal'}$ being optimal for $f$.
  }
\etheorem
\bproof
  If the algorithm lands on a minimizer of $f$ in a finite number of iterations, then there is nothing left to prove.  Thus, for the remainder of the proof, let us suppose that this does not occur.
  
  Our next aim is to show that, with the update \eqref{eq.delta_update}, one finds $\{\delta_k\} \searrow 0$.  To derive a contradiction, suppose that there exists $\delta \in (0,\delta_1]$ such that $\delta_k = \delta$ for all sufficient large $k \in \N{}$.  For $p \in (0,1)$, let $\Kcal \subseteq \N{}$ be the infinite index set for which Theorem~\ref{th.self-correction} and Corollary~\ref{cor.self-correction} guarantee the existence of $(\kappa,\sigma,\mu) \in \R{}_{>0} \times \R{}_{>0} \times \R{}_{>0}$ and $(\nu,\xi) \in \R{}_{>0} \times \R{}_{>0}$ such that \eqref{eq.good_bounds_H} and \eqref{eq.good_bounds} hold, and, by Theorem~\ref{th.svano}, such that \eqref{eq.svano_limits} holds.  By the optimality of $(\omega_k,\gamma_k)$ with respect to \eqref{prob.dual} for all $k \in \R{n}$, it follows that, for all $k \in \N{}$, there exists $t_k \in \R{n}$ such that
  \bequation\label{eq.overton2}
    s_k = -W_k(G_k\omega_k + \gamma_k) = \delta_kt_k,\ \ \|t_k\| \leq 1,\ \ \text{and}\ \ t_k^T\gamma_k = \|\gamma_k\|_*.
  \eequation
  (This fact is derived as \eqref{def.pk} in Appendix~\ref{app.primal-dual}.)  With $\pi_k$ representing the angle between $t_k$ and $\gamma_k$, the last equation in \eqref{eq.overton2} shows that there exists $c \in \R{}_{>0}$ such that
  \bequation\label{eq.reza}
    \cos(\pi_k) = t_k^T\gamma_k/(\|t_k\|_2\|\gamma_k\|_2) = \|\gamma_k\|_*/(\|t_k\|_2\|\gamma_k\|_2) \geq c/\|t_k\|_2
  \eequation
  for all $k \in \N{}$.  On the other hand, by the first equation in \eqref{eq.overton2}, the second limit in~\eqref{eq.svano_limits}, and the supposition that $\delta_k = \delta$ for all sufficiently large $k \in \N{}$, one finds that $\|t_k\|_2 \searrow 0$ over $k \in \Kcal$.  This limit and \eqref{eq.reza} imply that $\cos(\pi_k) \nearrow \infty$ over $k \in \Kcal$, a contradiction.  Hence, we may conclude that $\{\delta_k\} \searrow 0$.
  
  Since $\{\delta_k\} \searrow 0$, by construction of the update \eqref{eq.delta_update}, it follows that there exists an infinite index set~$\Kcal'$ such that \eqref{eq.delta_update} yields $\delta_{k+1} \gets \tau\delta_k$ for all $k \in \Kcal'$.  By \eqref{eq.delta_update}, this implies that \eqref{eq.almost_there} holds.  This limit, along with the fact that the subgradients $\{g_{k,j}\}$ used in~\ref{alg.bundle} are evaluated at points in $\Bmbb(x_k,\delta_k)$ for all $k \in \N{}$, implies by Theorem~\ref{th.stationary} that any limit point of $\{x_k\}_{k\in\Kcal'}$ is optimal.
\eproof

In \eqref{eq.delta_update}, one may replace the last term in the $\max$ with $\upsilon_3 \|G_k\omega_k\|_{W_k}$, which uses the $W_k$-norm rather than the Euclidean norm.  One can show with this update, as in the proof above, that $\{\delta_k\} \searrow 0$ and \eqref{eq.almost_there} holds for $\Kcal' := \{k \in \N{} : \delta_{k+1} < \delta_k\}$.  Indeed, with this modified update, a similar argument as in the proof above shows that $\{\delta_k\} \searrow 0$, implying that $\{\|G_k\omega_k\|_{W_k}\}_{k\in\Kcal'} \to 0$.  One then finds from \eqref{eq.overton2} that
\bequation\label{eq.referee}
  \delta_k \gamma_k^Tt_k = -\gamma_k^TW_kG_k\omega_k - \|\gamma_k\|_{W_k}^2 \implies \|\gamma_k\|_* = - \(\frac{\cos (\tilde\pi_k) \|\gamma_k\|_{W_k}\|G_k\omega_k\|_{W_k}}{\delta_k}\) - \frac{\|\gamma_k\|_{W_k}^2}{\delta_k},
\eequation
where $\cos (\tilde\pi_k) := \gamma_k^TW_kG_k\omega_k/(\|\gamma_k\|_{W_k}\|G_k\omega_k\|_{W_k})$.  By nonnegativity of norms and since $\delta_k > 0$ for all $k \in \N{}$, \eqref{eq.referee} shows that $-\cos (\tilde\pi_k) \geq 0$ for all $k \in \N{}$.  Moreover, for $k \in \Kcal'$, one finds from \eqref{eq.referee} and~\eqref{eq.delta_update} (again, with the last term in the max replaced by $\upsilon_3\|G_k\omega_k\|_{W_k}$) that
\bequation\label{eq.referee2}
  \|\gamma_k\|_* \leq -\(\frac{\cos(\tilde\pi_k)\|\gamma_k\|_{W_k}}{\upsilon_3}\) - \frac{\|\gamma_k\|_{W_k}^2}{\delta_k}.
\eequation
On the other hand, since $f$ is bounded below, $\{f(x_k)\}$ is monotonically decreasing, inequality \eqref{eq.sufficient_reduction} holds for all $k \in \N{}$, and $\{\|G_k\omega_k\|_{W_k}\}_{k\in\Kcal'} \to 0$, it follows that $\{\|\gamma_k\|_{W_k}\}_{k\in\Kcal'} \to 0$, which with \eqref{eq.referee2} implies that $\{\gamma_k\}_{k\in\Kcal'} \to 0$.  Combining this with $\{\|G_k\omega_k + \gamma_k\|_2\}_{k\in\Kcal'} \to 0$ shows that \eqref{eq.almost_there} holds.

%************
% Subsection
%************
\subsection{Bundle and Gradient Sampling Methods for Nonconvex Minimization}\label{subsec.nonconvex}

In this subsection, we describe two other methods that adhere to the~\ref{alg.svano} Framework that are designed for minimizing nonconvex and nonsmooth objective functions.  We do not present complete details and analyses for these methods, but direct the reader to articles for further information.

Building on \svanobundlestep{} in~\S\ref{subsec.bundle}, one can devise a bundle trust region method for minimizing nonconvex~$f$.  The ideas for doing so that we discuss here can be found in \S3, and in particular in ``Inner Iteration~(3.3),'' of \cite{SchrZowe92}.  Firstly, observe that when $f$ is nonconvex, one should not employ $b_{k,j}$ as defined in \eqref{def.b} since the mapping $x \mapsto f_{k,j} + g_{k,j}^T(x - x_{k,j})$ is no longer guaranteed to be an affine underestimator of $x \mapsto f(x)$.  Instead, one can ``downshift'' the cutting planes for the subproblem, in particular by defining the scalars instead to be
\bequation\label{def.b_downshift}
    b_{k,j} = \min\{f(x_k) - r\|x_k - x_{k,j}\|^2,f_{k,j} + g_{k,j}^T(x_k - x_{k,j})\}, \ \ \text{where}\ \ r \in \R{}_{>0}.
\eequation
This ensures that $f(x_k) \geq l_{k,m}(x_k)$, as required in the latter part of Lemma~\ref{lem.primal_dual}, and makes it less likely that the pair $(f_{k,j},g_{k,j})$ will influence the subproblem solution when $x_{k,j}$ is farther from $x_k$.  Two other modifications should also be made in the step computation procedure, one required for theoretical reasons and one due to practical considerations.
\bitemize
  \item For theoretical reasons, the algorithms proposed by \cite{SchrZowe92} sometimes involve the use of a ``Null Step'' in which a step of norm zero is accepted while (depending on certain verifiable conditions) either the local model of the objective is improved through the addition of a (downshifted) cutting plane or the trust region radius is reduced.  (In the context of~\ref{alg.svano}, one should \emph{skip} the update of the inverse Hessian approximation by setting $W_{k+1} \gets W_k$ when a null step occurs; the update should only be employed with a ``Serious Step".)  When $f$ is convex and the current iterate is not a minimizer, one can guarantee that after a series of null steps the algorithm will eventually compute a serious step that makes sufficient progress toward a minimizer.  However, when $f$ is nonconvex, one cannot always perform null steps and guarantee that a serious step that makes sufficient progress will ultimately be produced.  Instead, under certain verifiable conditions (indicating that one cannot guarantee sufficient progress either by reducing the trust region radius or adding a (downshifted) cutting plane), the proposed method falls back to a line search along a computed step.  (See Step~(2)(c) in ``Inner Iteration~(3.3)'' in \cite{SchrZowe92}.)  Such a fallback is also needed in the context of \ref{alg.svano}.
  \item For practical reasons (and to guarantee that a verifiable termination condition will ultimately be satisfied), one should not expend additional computational effort to produce a serious step if one has reasons to believe that the current iterate is approximately stationary.  In particular, in the context of \ref{alg.svano}, if $(\omega_k,\gamma_k,s_k)$ yields $\max\{\upsilon_1\|G_k\omega_k + \gamma_k\|_2,\upsilon_2\|s_k\|_2,\upsilon_3\|G_k\omega_k\|_2\} \leq \delta_k$, then the step computation should be terminated with a null step, one should \emph{skip} \eqref{step.W} by setting $W_{k+1} \gets W_k$, and the trust region radius should be decreased as in~\eqref{eq.delta_update}.  (See Step~(1) in ``Inner Iteration'' (3.3) in \cite{SchrZowe92}.)
\eitemize
The algorithm by \cite{SchrZowe92} also involves other features that attempt to avoid line searches, although one does not need to employ these features to maintain their convergence guarantees.  (Indeed, in our implementation used to obtain the results in \S\ref{sec.numerical}, we ran a weak Wolfe line search after each successful step.  This improved the performance of the algorithm, which we believe was a consequence of it computing better steps and pairs for the BFGS updating strategy.)  Critically important for the analysis in \cite{SchrZowe92} is an assumption that $f$, in addition to being locally Lipschitz, is also \emph{weakly (lower) semismooth}; see, e.g., \cite{Miff77}.  It is only with this assumption that one can guarantee that each line search will terminate finitely.

One can also devise gradient sampling strategies (see \cite{BurkLewiOver05} and \cite{Kiwi07}) that adhere to \ref{alg.svano}.  In particular, this is easily done with modifications to the method proposed by \cite{CurtQue15}.  This approach already employs inverse Hessian approximations computed by BFGS updates, although in a more restrictive manner in that an update is skipped if the curvature of the Hessian approximation along a computed step is not sufficiently large relative to the norm of the step and the computed stepsize is not sufficiently large.  (See Algorithm~5 in \cite{CurtQue15}.)  Within the \ref{alg.svano} Framework, this curvature condition does not need to be checked since the self-correcting properties of the updates will ensure that the inverse Hessian approximations have eigenvalues that are uniformly bounded below and above in a sufficient number of iterations.  (One should still skip an update if a computed stepsize is below a threshold, although, as shown in the analysis in \cite{CurtQue15}, the stepsize will be sufficiently large infinitely often if the algorithm does not terminate finitely.)  As for the bundle method described in the previous paragraph, one needs to assume that the objective function is weakly (lower) semismooth in order to guarantee that the line searches terminate finitely.

%*********
% Section
%*********
\section{Numerical Experiments}\label{sec.numerical}

We implemented \ref{alg.svano} in C++.  The software includes three algorithms: \bfgs, which follows the strategy in~\S\ref{subsec.bfgs} (with no trust region constraint for the subproblems); \svanobundle, which follows the strategy in~\S\ref{subsec.nonconvex} that builds on the method from~\S\ref{subsec.bundle}; and \svanogs, which, as described in \S\ref{subsec.nonconvex}, follows the adaptive gradient sampling strategy from \cite{CurtQue15} with the addition of a trust region constraint.  For solving the arising subproblems, the code has its own specialized active-set QP solver that borrows ideas from \cite{Kiwi86}, but also allows for the variable-metric induced by~$W_k$ (see \cite{CurtQue13}) and a trust region constraint with $\|\cdot\| = \|\cdot\|_\infty$, meaning $\|\cdot\|_* = \|\cdot\|_1$.

Our implementation of \ref{alg.svano} allows for various choices of $\Hbar$ and $y_k$ for all $k \in \N{}$.  For example, as often proposed for BFGS methods for smooth optimization (e.g., see~\cite{NoceWrig06}), one might choose the former as a multiple of the identity where the multiplying factor is determined by a \cite{BarzBorw88} ``two-point stepsize'' strategy, projected onto $[\eta,\theta]$, after the first accepted step.  One could even update it with each iteration, as long as the factor is projected onto $[\eta,\theta]$ for all $k \in \N{}$.  Another strategy would be to initialize $H_1 \gets I$, employ \eqref{eq.H_update} through iteration $K \in \N{}$, then set $\Hbar \gets H_{K+1}$ for use in all subsequent iterations.  As for $y_k$, one could choose the displacement $g_{k+1,1} - g_{k,1}$ or $g_{k+1,1} - G_k\omega_k$, where one should recall that $G_k\omega_k$ is a convex combination of subgradients.  However, for our experiments we simply set $W_1 \gets I$ and $y_k \gets g_{k+1,1} - g_{k,1}$ for all $k \in \N{}$ so that our comparison of the algorithms mentioned above would be based on common choices of these values.

Our code sets $\delta_1 \gets 1$, then sets the remaining values in the sequence $\{\delta_k\}$ according to \eqref{eq.delta_update} with inputs $\tau \gets 0.5$ and $\upsilon_1 = \upsilon_2 = \upsilon_3 \gets 1$.  These values are used for the trust region radii in \svanobundle{} and \svanogs{}, for the sampling radius in \svanogs{}, and in the termination conditions (see below) in all algorithms.  In addition, at the beginning of each iteration $k \in \N{}$ of \svanobundle{} and \svanogs{}, the code initializes $\{x_{k+1,j}\}$ with those points from $\{x_{k,j}\}$ that are within $\Bmbb(x_{k+1},\delta_{k+1})$.  The remaining inputs were $\alpha \gets 10^{-15}$, $\eta \gets 10^{-12}$, $\theta \gets 20$, and $r \gets 10^{-15}$.  The parameter $\theta$ had a large impact on performance; we discuss this more later in this section.

For test problems, we used the first ten from~\cite{HaarMietMaek04} with $n = 50$.  Pertinent information about the problems---namely, an indication of whether each problem is convex, the objective value at the initial point ($f(x_0)$), and the global minimum value of the objective ($f(x_*)$)---are given in Table~\ref{tab.problems}.  For further information, including the starting point for each problem; see~\cite{HaarMietMaek04}.

\btable[ht]
  \scriptsize
  \centering
  \texttt{
  \setlength{\tabcolsep}{12pt}
  \btabular{lcrr}
    \hline
    Name & Convex? & $f(x_0)$ & $f(x_*)$ \\
    \hline
    maxq               & Yes & 2500.0 &   0.0 \\
    mxhilb             & Yes &    4.5 &   0.0 \\
    chained lq         & Yes &   49.0 & -69.3 \\
    chained cb3 1      & Yes &  980.0 &  98.0 \\
    chained cb3 2      & Yes &  980.0 &  98.0 \\
    active faces       & No  &    3.9 &   0.0 \\
    brown function 2   & No  &   98.0 &   0.0 \\
    chained mifflin 2  & No  &  232.8 & -34.8 \\
    chained crescent 1 & No  &  292.3 &   0.0 \\
    chained crescent 2 & No  &  292.3 &   0.0 \\
    \hline
  \etabular
  }
  \caption{Test problem information for $n=50$.}
  \label{tab.problems}
\etable

Consistent with our theoretical analysis, the code terminates with a message of success when
\bequation\label{eq.termination}
  \|G_k\omega_k\|_2 \leq 10\delta_k\ \ \text{and}\ \ \delta_k \leq 10^{-4}\ \ \text{for some}\ \ k \in \N{}.
\eequation
The code terminates with a message of failure if the iteration limit of $10^4$ is reached or a computed stepsize is below $10^{-15}$, the latter playing the role of $\tilde\alpha_{\min}$ in Algorithm~\ref{alg.bfgs}.  The results obtained with these stopping conditions are shown in Table~\ref{tab.results}.  (It should be noted that since \svanogs{} is a randomized algorithm, its performance on a single problem can change from one run to the next, even from the same initial point.  In each of our experiments, we show the results of a single run for each problem.  In each case, we have observed that the results we provide are representative of the algorithm's average performance in general.)  In the table, an \texttt{Exit} of \texttt{Stationary} indicates a successful termination while that of \texttt{Iteration} or \texttt{Stepsize} indicates that the iteration or stepsize limit, respectively, was exceeded.  The values \texttt{$\delta_{\text{end}}$} and \texttt{$f(x_{\text{end}})$} indicate the final elements of the sequences $\{\delta_k\}$ and $\{f(x_k)\}$, and the counters \texttt{\#iter}, \texttt{\#func}, \texttt{\#grad}, and \texttt{\#subs} indicate the number of iterations, function evaluations, gradient evaluations, and subproblems solved.

\btable[ht]
  \scriptsize
  \centering
  \texttt{
  \setlength{\tabcolsep}{7pt}
  \btabular{llrrrrrr}
    \hline
    \multicolumn{8}{c}{\bfgs{}} \\
    \hline
    Name & Exit & $\delta_{\text{end}}$ & $f(x_{\text{end}})$ & \#iter & \#func & \#grad & \#subs \\
    \hline
              maxq &  Stationary &  +9.77e-05 &  +2.22e-07 &     457 &    1005 &     465 &     458 \\ 
            mxhilb &    Stepsize &  +1.56e-03 &  +4.37e-02 &     134 &    2243 &     162 &     135 \\ 
        chained lq &    Stepsize &  +5.00e-02 &  -6.93e+01 &     173 &    3107 &     175 &     174 \\ 
     chained cb3 1 &    Stepsize &  +1.00e-01 &  +9.80e+01 &     314 &    7061 &     315 &     315 \\ 
     chained cb3 2 &    Stepsize &  +1.00e-01 &  +9.82e+01 &     136 &    2783 &     144 &     137 \\ 
      active faces &    Stepsize &  +2.50e-02 &  +1.11e-15 &      26 &     716 &      29 &      27 \\ 
  brown function 2 &    Stepsize &  +1.00e-01 &  +2.02e-08 &     308 &    5245 &     309 &     309 \\ 
 chained mifflin 2 &    Stepsize &  +5.00e-02 &  -3.48e+01 &    1009 &   22080 &    1034 &    1010 \\ 
chained crescent 1 &    Stepsize &  +1.00e-01 &  +1.84e-01 &     121 &    4045 &     124 &     122 \\ 
chained crescent 2 &    Stepsize &  +1.00e-01 &  +3.40e-03 &     291 &    5470 &     292 &     292 \\     \hline
    \hline
    \multicolumn{8}{c}{\svanobundle{}} \\
    \hline
    Name & Exit & $\delta_{\text{end}}$ & $f(x_{\text{end}})$ & \#iter & \#func & \#grad & \#subs \\
    \hline
              maxq &  Stationary &  +9.77e-05 &  +8.14e-07 &     260 &     470 &     653 &     315 \\ 
            mxhilb &  Stationary &  +9.77e-05 &  +5.92e-05 &      61 &     756 &     712 &     116 \\ 
        chained lq &  Stationary &  +9.77e-05 &  -6.93e+01 &      14 &      85 &      85 &      66 \\ 
     chained cb3 1 &  Stationary &  +9.77e-05 &  +9.80e+01 &      19 &     158 &     170 &     146 \\ 
     chained cb3 2 &  Stationary &  +9.77e-05 &  +9.80e+01 &      30 &      71 &      85 &      48 \\ 
      active faces &  Stationary &  +9.77e-05 &  +5.21e-05 &      16 &     160 &     153 &      18 \\ 
  brown function 2 &  Stationary &  +9.77e-05 &  +2.05e-07 &      17 &      86 &      66 &      39 \\ 
 chained mifflin 2 &  Stationary &  +9.77e-05 &  -3.48e+01 &      57 &     927 &     961 &     899 \\ 
chained crescent 1 &  Stationary &  +9.77e-05 &  +8.20e-06 &      30 &     113 &      81 &      35 \\ 
chained crescent 2 &  Stationary &  +9.77e-05 &  +2.50e-06 &      98 &     902 &     976 &     870 \\ 
    \hline
    \hline
    \multicolumn{8}{c}{\svanogs{}} \\
    \hline
    Name & Exit & $\delta_{\text{end}}$ & $f(x_{\text{end}})$ & \#iter & \#func & \#grad & \#subs \\
    \hline
              maxq &  Stationary &  +9.77e-05 &  +8.17e-07 &     225 &     658 &     351 &     226 \\ 
            mxhilb &  Stationary &  +9.77e-05 &  +1.15e-04 &      96 &    1191 &     701 &     101 \\ 
        chained lq &  Stationary &  +9.77e-05 &  -6.93e+01 &      17 &     473 &     238 &      18 \\ 
     chained cb3 1 &  Stationary &  +9.77e-05 &  +9.80e+01 &    4023 &   32550 &    4491 &    4069 \\ 
     chained cb3 2 &  Stationary &  +9.77e-05 &  +9.80e+01 &     137 &     419 &     165 &     138 \\ 
      active faces &  Stationary &  +9.77e-05 &  +6.44e-03 &      16 &     367 &     188 &      17 \\ 
  brown function 2 &  Stationary &  +9.77e-05 &  +3.77e-02 &      19 &     504 &     462 &      20 \\ 
 chained mifflin 2 &  Stationary &  +9.77e-05 &  -3.48e+01 &    2116 &   24228 &    2991 &    2135 \\ 
chained crescent 1 &  Stationary &  +9.77e-05 &  +2.30e-05 &      38 &     211 &      66 &      40 \\ 
chained crescent 2 &  Stationary &  +9.77e-05 &  +6.65e-03 &    1475 &   11284 &    1636 &    1479 \\ 
    \hline
  \etabular
  }
  \caption{Termination status, solution properties, and counter values when \bfgs{}, \svanobundle{}, and \svanogs{} were employed to solve the test problems stated in Table~\ref{tab.problems}.}
  \label{tab.results}
\etable

One can see that \svanobundle{} and \svanogs{} behave quite well in the sense that they terminate with success for all problems.  Which algorithm performs the best for a particular problem depends on the performance measure of interest.  In particular, \svanobundle{} often requires fewer function evaluations, but sometimes at the expense of more gradient evaluations and subproblem solves.

More striking in the results in Table~\ref{tab.results} is the fact that \bfgs{} only terminated with a message of success for one problem.  For the remaining problems, the code terminated due to a small stepsize (below $10^{-15}$).  This provides evidence for our discussion in~\S\ref{subsec.bfgs}, where we stated that the main issue with proving convergence guarantees for a classical BFGS approach is that one cannot be sure that the stepsize would remain sufficiently large.  That being said, the final objective values yielded by \bfgs{} show that this code did not always perform poorly in terms of the final objective value!  For many problems, the final value was close to optimal.  (To try to verify approximate stationarity in practice, one could employ auxiliary procedures; see, e.g., \S6.3 in \cite{LewiOver13}.  The effects of this are seen for a code with results in Table~\ref{tab.other_codes} later on.)

To illustrate the benefits of self-correction, we also ran the experiments with the same settings, except with $\theta \gets \infty$, a choice that is \emph{not valid} in terms of ensuring our convergence guarantees.  The results with these inputs are given in Table~\ref{tab.results_no_self_correction}.  Clearly, the performance is not as good.  The final objective values are often nearly optimal, but the code has a difficult time satisfying our termination criteria.  We conjecture that this behavior can be explained as follows.  Firstly, observe from Corollary~\ref{cor.self-correction}, namely, inequality~\eqref{eq.good_bounds}, that the objective decrease that is guaranteed through \eqref{eq.sufficient_reduction} is at least
\bequation\label{eq.daniel}
  \thalf \alpha (G_k\omega_k+\gamma_k)^TW_k(G_k\omega_k+\gamma_k) \geq \thalf \alpha (\nu/\xi) \|W_k(G_k\omega_k+\gamma_k)\|_2^2 = \thalf \alpha (\nu/\xi) \|s_k\|_2^2.
\eequation
It is through these relationships that Theorem~\ref{th.svano} ensures that there exists a subsequence of iterations over which $\|G_k\omega_k+\gamma_k\|_2$ and $\|s_k\|_2$ vanish.  However, as $\theta \nearrow \infty$, one finds through the analysis in Appendix~\ref{app.self-correcting}---specifically, \eqref{eq.zeta_upper_beta} through the end of the proof of Theorem~\ref{th.self-correction}---that the values in Theorem~\ref{th.self-correction} have $\kappa \searrow 0$, $\sigma \searrow 0$, and $\mu \nearrow \infty$.  This, in turn, implies that $\nu \searrow 0$ and $\xi \nearrow \infty$, meaning that $\nu/\xi \searrow 0$ in \eqref{eq.daniel}.  Consequently, small reductions in the objective (which is all that can be obtained near a minimizer) might be obtained with relatively large $\|G_k\omega_k + \gamma_k\|_2$ and $\|s_k\|_2$.  This means that numerous steps could be taken near a minimizer until the trust region radius is reduced (recall \eqref{eq.delta_update}), which in turn means that numerous steps could be taken near a minimizer until the termination condition \eqref{eq.termination} is satisfied.
%, which is to say that even the ``good'' inverse Hessian approximations might be very ill-conditioned.  Consequently, while the algorithm might produce steps with small norm as a stationary point is approached, the algorithm might not recognize that the iterates it computed are approximately stationary since the convex combination $G_k\omega_k$ might remain large in norm.  Comparatively speaking, this is less likely to occur when $W_k$ is well-conditioned, since in that case the norm $\|s_k\| = \|W_k(G_k\omega_k + \gamma_k\|$ being small more often corresponds to $\|G_k\omega_k\|$ being small.  (That said, the inverse Hessian approximation updates are not overly restricted; the code allows the self-correcting properties of BFGS to ensure this correspondence.)
Overall, this influence of $\theta$ on our numerical results shows that a practical benefit of our self-correcting framework is that it allows our code to enforce theoretically sound termination criteria.

In theory, a small value for $\eta$ might have a similar effect as a large value for $\theta$ as described in the previous paragraph.  However, we did not see the same effect in our experiments.  Indeed, the value that we used since it worked well in our experiments, namely, $\eta \gets 10^{-12}$, is relatively small.  One possible theoretical explanation for this is the fact that it enters as $\ln\eta$ in \eqref{eq.zeta_upper_beta}, whereas $\theta$ itself enters.

\btable[ht]
  \scriptsize
  \centering
  \texttt{
  \setlength{\tabcolsep}{7pt}
  \btabular{llrrrrrr}
    \hline
    \multicolumn{8}{c}{\bfgs{} (not enforcing \eqref{eq.bounds})} \\
    \hline
    Name & Exit & $\delta_{\text{end}}$ & $f(x_{\text{end}})$ & \#iter & \#func & \#grad & \#subs \\
    \hline
              maxq &  Stationary &  +9.77e-05 &  +2.35e-07 &     513 &    1392 &     611 &     514 \\ 
            mxhilb &   Iteration &  +1.56e-03 &  +9.31e-09 &   10000 &  478038 &   10145 &   10001 \\ 
        chained lq &    Stepsize &  +5.00e-02 &  -6.93e+01 &     343 &    2449 &     346 &     344 \\ 
     chained cb3 1 &    Stepsize &  +1.00e-01 &  +9.80e+01 &     500 &    2701 &     501 &     501 \\ 
     chained cb3 2 &    Stepsize &  +1.00e-01 &  +9.80e+01 &    2004 &   60334 &    2053 &    2005 \\ 
      active faces &    Stepsize &  +2.50e-02 &  +4.88e-15 &      37 &     160 &      40 &      38 \\ 
  brown function 2 &    Stepsize &  +1.00e-01 &  +2.82e-10 &     118 &     554 &     119 &     119 \\ 
 chained mifflin 2 &   Iteration &  +5.00e-02 &  -3.48e+01 &   10000 &  192578 &   11178 &   10001 \\ 
chained crescent 1 &    Stepsize &  +1.00e-01 &  +4.23e-11 &      72 &     754 &      99 &      73 \\ 
chained crescent 2 &    Stepsize &  +1.00e-01 &  +3.10e-14 &     546 &    2929 &     551 &     547 \\ 
    \hline
    \hline
    \multicolumn{8}{c}{\svanobundle{} (not enforcing \eqref{eq.bounds})} \\
    \hline
    Name & Exit & $\delta_{\text{end}}$ & $f(x_{\text{end}})$ & \#iter & \#func & \#grad & \#subs \\
    \hline
              maxq &  Stationary &  +9.77e-05 &  +5.11e-07 &     176 &     456 &     477 &     203 \\ 
            mxhilb &    Stepsize &  +3.91e-04 &  +8.96e-06 &      79 &     582 &     331 &     100 \\ 
        chained lq &  Stationary &  +9.77e-05 &  -6.93e+01 &      15 &     557 &     374 &     356 \\ 
     chained cb3 1 &    Stepsize &  +2.50e-02 &  +9.80e+01 &     205 &    4758 &    4867 &    4658 \\ 
     chained cb3 2 &   Iteration &  +1.25e-02 &  +9.80e+01 &   10000 &  897132 &  508499 &  498494 \\ 
      active faces &  Stationary &  +9.77e-05 &  +9.51e-05 &      20 &     465 &     383 &      24 \\ 
  brown function 2 &  Stationary &  +9.77e-05 &  +1.20e-09 &      17 &     435 &     369 &      74 \\ 
 chained mifflin 2 &    Stepsize &  +1.25e-02 &  -3.48e+01 &      75 &    1245 &    1116 &    1033 \\ 
chained crescent 1 &  Stationary &  +9.77e-05 &  +2.50e-09 &      72 &     468 &     383 &     284 \\ 
chained crescent 2 &    Stepsize &  +1.25e-02 &  +1.56e-03 &     144 &    4680 &    4215 &    4051 \\ 
    \hline
    \hline
    \multicolumn{8}{c}{\svanogs{} (not enforcing \eqref{eq.bounds})} \\
    \hline
    Name & Exit & $\delta_{\text{end}}$ & $f(x_{\text{end}})$ & \#iter & \#func & \#grad & \#subs \\
    \hline
              maxq &  Stationary &  +9.77e-05 &  +1.01e-06 &     141 &     439 &     242 &     142 \\ 
            mxhilb &  Stationary &  +9.77e-05 &  +6.28e-07 &     115 &     768 &     409 &     120 \\ 
        chained lq &  Stationary &  +9.77e-05 &  -6.93e+01 &      17 &     755 &     491 &      18 \\ 
     chained cb3 1 &   Iteration &  +1.25e-02 &  +9.88e+01 &   10000 &  411205 &    2720 &   10002 \\ 
     chained cb3 2 &   Iteration &  +3.13e-03 &  +9.80e+01 &   10000 &  504074 &     147 &   10001 \\ 
      active faces &  Stationary &  +9.77e-05 &  +6.44e-03 &      16 &     440 &     353 &      18 \\ 
  brown function 2 &  Stationary &  +9.77e-05 &  +3.77e-02 &      19 &     551 &     461 &      20 \\ 
 chained mifflin 2 &   Iteration &  +6.25e-03 &  -3.47e+01 &   10000 &  459942 &    4333 &   10047 \\ 
chained crescent 1 &  Stationary &  +9.77e-05 &  +1.06e-05 &      61 &     171 &      85 &      62 \\ 
chained crescent 2 &   Iteration &  +7.81e-04 &  +2.29e-02 &   10000 &  495522 &     758 &   10032 \\ 
    \hline
  \etabular
  }
  \caption{Termination status, solution properties, and counter values when \bfgs{}, \svanobundle{}, and \svanogs{} were employed to solve the test problems stated in Table~\ref{tab.problems}.  Unlike for Table~\ref{tab.results}, these results were obtained with $\theta \gets \infty$, meaning that the latter bound in \eqref{eq.bounds} is not enforced.}
  \label{tab.results_no_self_correction}
\etable

For reference, we provide in Table~\ref{tab.other_codes} the results when solving the problems with \texttt{LMBM}\footnote{http://napsu.karmitsa.fi/lmbm/} (written in Fortran) and \texttt{HANSO}\footnote{https://cs.nyu.edu/overton/software/hanso/} (written in Matlab), both using their default settings.  It is difficult to compare the performance of these codes with our methods since the termination conditions for all codes are different.  Indeed, while our methods \emph{only} terminate with a message of success if \eqref{eq.termination} is satisfied, \texttt{LMBM} and \texttt{HANSO} terminate due to other sets of conditions.  (See the caption of Table~\ref{tab.other_codes}.)  Variants of SVANO might benefit from tailored termination conditions depending on the type of algorithm (e.g., a bundle versus a gradient sampling method) and class of problem being solved, but in the interest of having consistent experiments based on our general theoretical results, we have required \eqref{eq.termination}.  (It is worth mentioning that for the problems that \svanogs{} struggled to solve---\texttt{chained~cb3~1}, \texttt{chained~mifflin~2}, and \texttt{chained~crescent~2}---both \texttt{LMBM} and \texttt{HANSO} terminated with \texttt{Exit} not equal to 1.)

\btable[ht]
  \scriptsize
  \centering
  \texttt{
  \setlength{\tabcolsep}{7pt}
  \btabular{l|crrr|crrr}
    \hline
    & \multicolumn{4}{c|}{\texttt{LMBM}} & \multicolumn{4}{c}{\texttt{HANSO}} \\
    \hline
    Name & Exit & $f(x_{\text{end}})$ & \#iter & \#func & Exit & $f(x_{\text{end}})$ & \#iter & \#func \\
    \hline
              maxq & 1 & +4.97e-06 & 460 &  501 & 1 & +1.18e-08 & 494 & 1002 \\
            mxhilb & 3 & +1.25e-06 & 377 & 1385 & 2 & +1.92e-12 & 473 & 1091 \\
        chained lq & 3 & -6.93e+01 & 219 & 1233 & 1 & -6.93e+01 & 182 & 811 \\
     chained cb3 1 & 3 & +9.80e+01 & 167 &  750 & 3 & +9.80e+01 & 149 & 968 \\
     chained cb3 2 & 3 & +9.81e+01 &  20 &   61 & 3 & +9.80e+01 &  83 & 208 \\
      active faces & 2 & +7.40e-11 &  96 &   97 & 1 & +2.35e-05 &  11 &  27 \\
  brown function 2 & 2 & +1.36e-08 & 348 & 3195 & 1 & +1.24e-04 &  24 & 102 \\
 chained mifflin 2 & 3 & -3.48e+01 & 258 & 1890 & 3 & -3.48e+01 & 962 & 2858 \\
chained crescent 1 & 3 & +2.65e-09 & 106 &  294 & 1 & +9.45e-06 &  21 & 51 \\
chained crescent 2 & 2 & +6.12e-05 & 488 & 4351 & 3 & +7.40e-07 &  92 & 453 \\
    \hline
  \etabular
  }
  \caption{Termination status, final objective value, and counter values when \texttt{LMBM} and \texttt{HANSO} are employed to solve the test problems stated in Table~\ref{tab.problems}.  For \texttt{LMBM}, an \texttt{Exit} of 1 means ``the problem has been solved with desired accuracy,'' while an \text{Exit} of 2 or 3 means that changes in the objective were sufficiently small.  (See \texttt{IOUT(3)} in \texttt{LMBM}'s documentation.)  For \texttt{HANSO}, an \texttt{Exit} of 1 means ``norm of smallest vector in convex hull of gradients below tolerance,'' of 2 means a direction of ascent was computed, and of 3 means that the line search failed.}
  \label{tab.other_codes}
\etable

%*********
% Section
%*********
\section{Conclusion}\label{sec.conclusion}

We have proposed a framework for solving nonsmooth optimization problems.  Its distinguishing characteristic is that it maintains and benefits from the self-correcting properties of BFGS updating of the generated sequence of inverse Hessian approximations.  In particular, it benefits theoretically in that global convergence guarantees can be established, and it benefits in practice in that instances of the framework are effectively able to determine when iterates are nearly stationary for the objective.

Our discussions and analysis have been presented under Assumption~\ref{ass.f}.  One might also be interested in situations when~$f$ can be unbounded below and/or when it is extended-real-valued, i.e., when $f : \R{n} \to (\R{} \cup \{-\infty,\infty\})$.  We claim that the proposed framework, which ensures monotonic decrease in~$f$, is also viable in such cases, at least as long as one has access to an initial iterate $x_1$ in the effective domain of $f$, i.e., $x_1 \in \dom(f) := \{x \in \R{n} : f(x) < \infty\}$.  If $f$ is unbounded below and an iterate sequence $\{x_k\}$ is generated such that $\{f(x_k)\} \searrow -\infty$, then there is nothing else that one should ask from the proposed framework.  Hence, for simplicity, our Assumption~\ref{ass.f} precluded this case by ensuring that any such sequence~$\{f(x_k)\}$ is bounded below.  As for cases when $f$ is extended-real-valued, we claim that if any stationary point for~$f$ lies in the interior of the effective domain $\dom(f)$, then, with slight modifications of the proposed framework---e.g., to handle points encountered outside $\dom(f)$---our analysis for our framework follows in essentially the same manner as under Assumption~\ref{ass.f}.  

It is worthwhile to point out that we have not discussed limited-memory BFGS, even though using limited memory ideas is another alternative for ensuring that the inverse Hessian approximations have eigenvalues that are uniformly bounded below (away from zero) and above.  The primary reason for this omission is that we have observed that limited memory BFGS techniques do not typically perform as well as a full memory approach in the context of nonsmooth optimization.

%*********
% Section
%*********
\section*{Acknowledgements}\label{sec.acknowledgement}

This material is based upon work supported by the U.S.~Department of Energy, Office of Science, Applied Mathematics, Early Career Research Program under Award Number DE--SC0010615 and by the U.S.~National Science Foundation, Division of Mathematical Sciences, Computational Mathematics Program under Award Numbers DMS--1016291 and DMS--1319356.

The authors would like to thank Andreas W\"achter for hosting and guiding the third author during that author's visit to Northwestern University in the summer of 2017 while that author was implementing the quadratic optimization solver for our software.  They would also like to thank Michael L.~Overton for providing numerous valuable comments that helped to improve a draft of the paper.  Last, but not least, the authors would like to thank the anonymous referees whose very thoughtful reports and interesting exchanges with us led to nice improvements to the paper.

%\newpage

%**********
% Appendix
%**********
\appendix
\numberwithin{equation}{section}

%*********
% Section
%*********
\section{Primal and Dual Subproblems}\label{app.primal-dual}

In this appendix, we show that the dual of~\eqref{prob.primal} is~\eqref{prob.dual}, how the solution of \eqref{prob.primal} can be recovered from that of \eqref{prob.dual}, and that Lemma~\ref{lem.primal_dual} holds true.

As previously mentioned in~\S\ref{sec.step_computations}, the primal problem~\eqref{prob.primal} is equivalent to \eqref{prob.primal_smooth}.  A Lagrangian for this problem, call it $L : \R{n} \times \R{} \times \R{m} \to \R{}$, is given by
\bequationNN
  L(x,z,\omega) = z + \thalf (x - x_k)^TH_k(x - x_k) + \sum_{j=1}^m \omega_j(f_{k,j} + g_{k,j}^T(x - x_{k,j}) - z),
\eequationNN
with which we can write the dual problem for \eqref{prob.primal_smooth} (see \cite{Bert09}) as
\bequationNN
  \sup_{\omega \in \R{m}_+}\ \inf_{(x,z) \in \Xcal_k \times \R{}}\ L(x,z,\omega).
\eequationNN
Differentiating $L$ with respect to $z$, one finds that the ``inner'' infimum is attained only if $\mathds{1}^T\omega = 1$, from which it follows that the dual is equivalent to
\bequation\label{prob.dual_supinf}
  \sup_{\omega \in \R{m}_+}\ \(\inf_{x \in \Xcal_k}\ \(\thalf (x - x_k)^TH_k(x - x_k) + \sum_{j=1}^m \omega_j(f_{k,j} + g_{k,j}^T(x - x_{k,j}))\)\)\ \ \st\ \ \mathds{1}^T\omega = 1.
\eequation
Defining the characteristic $\chi_{\Xcal_k} : \R{n} \to \R{} \cup \{\infty\}$ as one that evaluates as 0 for $x \in \Xcal_k$ and $\infty$ otherwise, the inner infimum problem can equivalently be written as
\bequation\label{prob.inner_inf}
  \inf_{x \in \R{n}}\ \Lbar(x) + \chi_{\Xcal_k}(x),
\eequation
where we define the quadratic function $\Lbar : \R{n} \to \R{}$ by
\bequationNN
  \baligned
  \Lbar(x) 
  &= \thalf (x - x_k)^TH_k(x - x_k) + \sum_{j=1}^m \omega_j(f_{k,j} + g_{k,j}^T(x - x_{k,j})) \\
  &= \thalf x^TH_kx + x^T\(-H_kx_k + \sum_{j=1}^m \omega_j g_{k,j}\) + \thalf x_k^TH_kx_k + \sum_{j=1}^m \omega_j(f_{k,j} - g_{k,j}^Tx_{k,j}).
  \ealigned
\eequationNN
The conjugate of $\Lbar$, namely $\Lbar^\star : \R{n} \to \R{}$, is given by\footnote{Recall that for $A \in \R{n \times n}$, $b \in \R{n}$, and $c \in \R{}$ with $A \succ 0$, the conjugate of $\Lbar : \R{n} \to \R{}$ defined by $\Lbar(x) = \thalf x^TAx + b^Tx + c$ is given by $\phi^\star(y) = \thalf(y-b)^TA^{-1}(y-b) - c$.  For example, see~\cite{Bert09}.}
\bequationNN
  \baligned
  \Lbar^\star(y)
    =&\ \thalf\(y + \(H_kx_k - \sum_{j=1}^m \omega_jg_{k,j}\)\)^TW_k\(y + \(H_kx_k - \sum_{j=1}^m \omega_jg_{k,j}\)\) \\
     &\ -\thalf x_k^TH_kx_k - \sum_{j=1}^m \omega_j(f_{k,j} - g_{k,j}^Tx_{k,j}) \\
    =&\ \thalf\(y - \sum_{j=1}^m \omega_jg_{k,j}\)^TW_k\(y - \sum_{j=1}^m \omega_jg_{k,j}\) + x_k^Ty - \sum_{j=1}^m \omega_j(f_{k,j} + g_{k,j}^T(x_k - x_{k,j})).
  \ealigned
\eequationNN
In addition, the conjugate of $\chi_{\Xcal_k}$, namely $(\chi_{\Xcal_k})^\star : \R{n} \to \R{}$, is given by
\bequationNN
  (\chi_{\Xcal_k})^\star(y) = \sup_{x\in\R{n}} (y^Tx - \chi_{\Xcal_k}(x)) = \sup_{x\in\Xcal_k} y^Tx = \sup_{\|x - x_k\| \leq \delta_k} y^Tx.
\eequationNN
If $\delta_k = \infty$, then $(\chi_{\Xcal_k})^\star(y) = \infty$ for all nonzero $y \in \R{n}$. Otherwise, defining the vector $s := (x - x_k)/\delta_k$ so that $x = \delta_ks + x_k$, the above implies that
\bequationNN
  (\chi_{\Xcal_k})^\star(y) = \sup_{\|x - x_k\| \leq \delta_k} y^Tx = \sup_{\|s\|\leq 1} y^T(\delta_k s + x_k) = x_k^Ty + \delta_k\|y\|_\star.
\eequationNN
In either case, since the intersection of the relative interiors of the effective domains of~$\Lbar$ and~$\chi_k$ is nonempty, Fenchel duality implies the strong duality relationship
\bequationNN
  \baligned
  \inf_{x\in\R{n}}\ \Lbar(x) + \chi_{\Xcal_k}(x) = \sup_{y\in\R{n}} & -\Lbar^\star(y) - (\chi_{\Xcal_k})^\star(-y) \\
  = \sup_{y\in\R{n}} & -\thalf\(y - \sum_{j=1}^m \omega_jg_{k,j}\)^TW_k\(y - \sum_{j=1}^m \omega_jg_{k,j}\) \\
  & + \sum_{j=1}^m \omega_j(f_{k,j} + g_{k,j}^T(x_k - x_{k,j})) - \delta_k\|y\|_\star,
  \ealigned
\eequationNN
where, for the case $\delta_k = \infty$, we set $y = 0$ and interpret $\delta_k\|y\|_\star$ as zero. Going back to \eqref{prob.dual_supinf}, we now deduce that this problem is equivalent to
\bequationNN
  \baligned
  \sup_{(\omega,y) \in \R{m}_+ \times \R{n}}
  &\ -\thalf\(y - \sum_{j=1}^m \omega_jg_{k,j}\)^TW_k\(y - \sum_{j=1}^m \omega_jg_{k,j}\) + \sum_{j=1}^m \omega_j(f_{k,j} + g_{k,j}^T(x_k - x_{k,j})) - \delta_k\|y\|_\star \\
  \st &\ \mathds{1}^T\omega = 1.
  \ealigned
\eequationNN
Letting $\gamma = -y$ and observing \eqref{def.b}, this leads to \eqref{prob.dual}, as desired.

\bproof[Proof of Lemma~\ref{lem.primal_dual}]
Let us show that with $(\omega_k,\gamma_k)$ solving \eqref{prob.dual} and $s_k$ defined in \eqref{eq.s}, the point~$x_{k+1}$ in \eqref{eq.x_update} solves \eqref{prob.primal}.  Firstly, optimality of $(\omega_k,\gamma_k)$ implies that, with the optimal $z_k \in \R{}$ for \eqref{prob.primal_smooth},
\bsubequations
  \begin{align}
    0 &=   -G_k^TW_k(G_k\omega_k + \gamma_k) + b_k - z_k\mathds{1} \label{eq.dual.kkt_omega} \\ \text{and}\ \ 
    0 &\in -W_k(G_k\omega_k + \gamma_k) - \delta_k \partial \|\gamma_k\|_*.
  \end{align}
\esubequations
Using the fact that for any $\gamma \in \R{n}$ one has
\bequationNN
  \partial \|\gamma\|_* = \{t \in \R{n} : \|t\| \leq 1\ \text{and}\ t^T\gamma = \|\gamma\|_*\},
\eequationNN
it follows that there exists a vector $t_k \in \R{n}$ such that
\bequation\label{def.pk}
  0 = W_k(G_k\omega_k + \gamma_k) + \delta_k t_k,\ \ \|t_k\| \leq 1,\ \ \text{and}\ \ t_k^T \gamma_k = \|\gamma_k\|_*.
\eequation
Hence, evaluating the dual objective function at $(\omega_k,\gamma_k)$, one obtains
\bequationNN
  \baligned
   &\ -\thalf(G_k\omega_k+\gamma_k)^TW_k(G_k\omega_k+\gamma_k) + b_k^T\omega_k - \delta_k\|\gamma_k\|_* \\
  =&\ -\thalf (G_k \omega_k+\gamma_k)^T W_k (G_k\omega_k+\gamma_k) + b_k^T\omega_k - \delta_k t_k^T \gamma_k \\
  =&\ -\thalf (G_k \omega_k+\gamma_k)^T W_k (G_k\omega_k+\gamma_k) + b_k^T\omega_k + \gamma_k^TW_k(G_k\omega_k+\gamma_k) \\
  =&\ \thalf (G_k \omega_k+\gamma_k)^T W_k (G_k\omega_k+\gamma_k) + b_k^T\omega_k - \omega_k^TG_k^T W_k(G_k\omega_k+\gamma_k).
  \ealigned
\eequationNN
By duality theory, our desired conclusion follows as long as $x_{k+1}$ is feasible for~\eqref{prob.primal} and yields an objective value equal to this dual objective value.  To see that this is the case, first notice that, by \eqref{def.pk},
\bequationNN
  \|x_{k+1} - x_k\| = \|s_k\| = \|W_k(G_k\omega_k + \gamma_k)\| = \delta_k\|t_k\| \leq \delta_k.
\eequationNN
Secondly, observe that the primal objective value at $x_{k+1}$ is
\bequationNN
  \baligned
  q_{k,m}(x_{k+1})
  &= \thalf(x_{k+1} - x_k)^TH_k(x_{k+1} - x_k) + z_k \\
  &= \thalf(G_k\omega_k + \gamma_k)^TW_k(G_k\omega_k + \gamma_k) + b_k^T\omega_k - \omega_k^TG_k^TW_k(G_k\omega_k + \gamma_k),
  \ealigned
\eequationNN
where the second equation follows from the definition of $s_k$ in \eqref{eq.s}, the result of multiplying \eqref{eq.dual.kkt_omega} on the left by $\omega_k^T$, and the fact that $\omega_k^T\mathds{1} = 1$.

All that remains is to prove that inequality \eqref{eq.shim} holds under the assumption that 
$f(x_k) \geq l_{k,m}(x_k)$. 
Since $x_k$ is feasible for \eqref{prob.primal} yielding an objective value of $q_{k,m}(x_k) = l_{k,m}(x_k) \leq f(x_k)$,
it follows that
\bequationNN
  0 \leq f(x_k) - q_{k,m}(x_{k+1}) = f(x_k) - l_{k,m}(x_{k+1}) - \thalf(x_{k+1} - x_k)^TH_k(x_{k+1} - x_k),
\eequationNN
from which it follows that
\bequationNN
  f(x_k) - l_{k,m}(x_{k+1}) \geq \thalf(x_{k+1} - x_k)^TH_k(x_{k+1} - x_k) = \thalf (G_k \omega_k+\gamma_k)^T W_k (G_k\omega_k+\gamma_k),
\eequationNN
yielding \eqref{eq.shim}, as desired.
\eproof

%*********
% Section
%*********
\section{Geometric Properties of BFGS Updating}\label{sec.geometric}

The update~\eqref{eq.H_update} performs a \emph{projection} to eliminate certain curvature dictated by the Hessian approximation $H_k$, as well as a corresponding \emph{correction} that replaces this curvature for the new approximation~$H_{k+1}$.  The details of this projection and associated correction in the Hessian approximation can be seen in the following manner.\footnote{The presentation in this appendix is based on notes by James V.~Burke, Adrian S.~Lewis, and Michael L.~Overton, which were shared with the first author by Michael L.~Overton.}  For the sake of generality, let $H \in \R{n\times n}$ be symmetric positive definite, i.e., $H \succ 0$.  An inner product based on $H$ (i.e., an ``$H$-inner product'') is
\bequationNN
  \inprod{s}{v}_H := s^THv\ \ \text{for all}\ \ (s,v) \in \R{n} \times \R{n}.
\eequationNN
Given a matrix $A \in \R{n \times n}$, let its ``$H$-adjoint'' with respect to this $H$-inner product, call it $A^* \in \R{n\times n}$, be a matrix that satisfies
\bequationNN
  \inprod{s}{Av}_H = \inprod{A^*s}{v}_H\ \ \text{for all}\ \ (s,v) \in \R{n} \times \R{n}.
\eequationNN
Since $H \succ 0$, it is easily verified that the unique $H$-adjoint of $A$ is $A^* = H^{-1}A^TH$.  One calls $A$ an ``$H$-orthogonal'' projection matrix if and only if it satisfies
\bequationNN
  A = A^2\ \ \text{(i.e., $A$ is idempotent)}\ \ \text{and}\ \ A = A^*\ \ \text{(i.e., $A$ is ``$H$-self-adjoint'')}.
\eequationNN
For example, given a nonzero vector $s \in \R{n}$, consider the matrices
\bequationNN
  P := \frac{ss^TH}{s^THs}\ \ \text{and}\ \ Q := I - P.
\eequationNN
Note that $P$ and $Q$ are $H$-orthogonal projection matrices; i.e., $P$ yields the $H$-orthogonal projection onto $\linspan(s)$ while $Q$ yields the $H$-orthogonal projection onto the subspace $H$-orthogonal to $\linspan(s)$.  That is, given $t \in \R{n}$ such that $\inprod{s}{t}_H = 0$ (i.e., $t$ lies in the subspace $H$-orthogonal to $\linspan(s)$), it follows that
\bequationNN
  \left\{ \baligned Ps &= s \\ Qs &= 0 \ealigned \right\}
  \ \ \text{while}\ \ 
  \left\{ \baligned Pt &= 0 \\ Qt &= t \ealigned \right\}.
\eequationNN

One may now interpret the updates yielded by \eqref{eq.H_update} in terms of sequences of projections and corrections.  Specifically, note that \eqref{eq.H_update} can be rewritten as
\bequationNN
  H_{k+1} \gets H_k^{(n-1)} + H_k^{(1)}\ \ \text{where}\ \ H_k^{(n-1)} := Q_k^TH_kQ_k\ \ \text{with}\ \ Q_k := \(I - \frac{s_ks_k^TH_k}{s_k^TH_ks_k}\)\ \ \text{and}\ \ H_k^{(1)} := \frac{v_kv_k^T}{s_k^Tv_k}.
\eequationNN
Based on the discussion above, $Q_k$ yields the $H_k$-orthogonal projection onto the subspace $H_k$-orthogonal to $\linspan(s_k)$.  Looking more closely, $H_k^{(n-1)}$ has rank $n-1$, remains positive definite on the subspace $H$-orthogonal to $\linspan(s_k)$, and
\bequation\label{eq.projected}
  H_k^{(n-1)}s_k = 0\ \ \text{while}\ \ H_k^{(n-1)}t = H_kt\ \ \text{if}\ \ \inprod{s_k}{t}_H = 0.
\eequation
On the other hand, the matrix $H_k^{(1)}$ can be written as
\bequationNN
  H_k^{(1)} = \frac{v_kv_k^T}{s_k^Tv_k} = \(\frac{\|v_k\|_2^2}{s_k^Tv_k}\)\(\frac{v_kv_k^T}{\|v_k\|_2^2}\),
\eequationNN
where, in light of the secant-like equation \eqref{eq.secant}, the leading scalar $\|v_k\|_2^2/s_k^Tv_k = \|v_k\|_2^2/v_k^TW_{k+1}v_k$ is the inverse of a Rayleigh quotient for $W_{k+1}$ defined by $v_k$.  Since the so-called ``curvature condition'' $s_k^Tv_k > 0$ holds by~\eqref{eq.bounds}, one finds that $s_k^TH_k^{(1)}s_k = s_k^Tv_k > 0$, so one finds that $H_k^{(1)}$ corrects the curvature along $\linspan(s_k)$ that, according to \eqref{eq.projected}, has been projected out of $H_k^{(n-1)}$.

%*********
% Section
%*********
\section{Self-Correcting Properties of BFGS Updating}\label{app.self-correcting}

The purpose of this appendix is to provide a proof of Theorem~\ref{th.self-correction}.  Firstly, observe that nonpositivity of the latter terms in \eqref{eq.gamma} follows since $s_k \neq 0$ and $H_k \succ 0$ imply that $\cos^2\phi_k \in (0,1]$ and $\iota_k \in \R{}_{>0}$ for all $k \in \N{}$, since $\ln(r) \leq 0$ for all $r \in (0,1]$, and since $1 - r + \ln(r) \leq 0$ for all $r \in \R{}_{>0}$.  Note also that
\bequation\label{eq.zeta}
  1 - r + \ln(r) \searrow -\infty\ \ \text{as}\ \ r \searrow 0\ \ \text{or}\ \ r \nearrow \infty.
\eequation
These facts are used explicitly in the following proof.

\bproof[Proof of Theorem~\ref{th.self-correction}]
  For all $k \in \N{}$, define
  \bequation\label{eq.zeta2}
    \zeta_k := -\ln(\cos^2\phi_k) - \(1 - \frac{\iota_k}{\cos^2\phi_k} + \ln\(\frac{\iota_k}{\cos^2\phi_k}\)\) \geq 0.
  \eequation
  Hence, from \eqref{eq.bounds} and \eqref{eq.gamma}, it follows for all $k \in \N{}$ that
  \bequationNN
    \psi(H_{k+1}) \leq \psi(H_k) + \theta - 1 - \ln \eta - \zeta_k,
  \eequationNN
  from which it follows for all $K \in \N{}$ that
  \bequationNN
    \psi(H_{K+1}) \leq \psi(H_1) + (\theta - 1 - \ln\eta)K - \sum_{k=1}^K \zeta_k.
  \eequationNN
  Since, for all $K \in \N{}$, one has $\psi(H_{K+1}) \in \R{}_{>0}$, this implies that
  \bequation\label{eq.zeta_upper}
    \frac{1}{K}\sum_{k=1}^K \zeta_k < \frac{1}{K}\psi(H_1) + (\theta - 1 - \ln\eta).
  \eequation
  Now, considering fixed $p \in (0,1)$ and $K \in \N{}$, let $J_{p,K}$ be the set of indices corresponding to the $\lceil pK \rceil$ smallest elements of~$\zeta_k$ for $k \in \{1,\dots,K\}$, and let $\zeta_{p,K}$ denote the largest element of $\{\zeta_k\}_{k \in J_{p,K}}$.  Then,
  \bequationNN
    \frac{1}{K}\sum_{k=1}^K \zeta_k \geq \frac{1}{K}\(\zeta_{p,K} + \sum_{k=1,k \notin J_{p,K}}^K \zeta_k\) \geq \(\frac{1}{K} + \frac{K - \lceil pK \rceil}{K}\)\zeta_{p,K} \geq (1-p)\zeta_{p,K},
  \eequationNN
  which along with \eqref{eq.zeta_upper} and the fact that $K \geq 1$ implies that, for all $k \in J_{p,K}$,
  \bequation\label{eq.zeta_upper_beta}
    \zeta_k \leq \zeta_{p,K} < \frac{1}{1-p}(\psi(H_1) + \theta - 1 - \ln\eta) =: c_0 \in \R{}_{>0}.
  \eequation
  Since the facts that $\cos^2\phi_k \in (0,1]$ and $\iota_k \in \R{}_{\geq0}$ for all $k \in \N{}$, \eqref{eq.zeta},
$1 - r + \ln(r) \leq 0$ for all $r \in \R{}_{>0}$, 
and \eqref{eq.zeta2} together imply that $\zeta_k \geq -\ln(\cos^2\phi_k)$ for all $k \in J_{p,K}$, it follows from \eqref{eq.zeta_upper_beta} that $-\ln(\cos^2\phi_k) < c_0$ for all $k \in J_{p,K}$, which means that $\cos \phi_k > e^{-c_0/2} =: c_1 \in \R{}_{>0}$ for all $k \in J_{p,K}$.  That is, observing \eqref{eq.cos}, the first inequality in \eqref{eq.good_bounds_H} holds for any constant $\kappa \in (0,c_1]$ and for all $k\in J_{p,K}$.  Now observe that \eqref{eq.zeta2}, the fact that $-\ln(\cos^2\phi_k) \geq 0$ for all $k \in \N{}$, and \eqref{eq.zeta_upper_beta} imply for all $k \in J_{p,K}$ that
  \bequationNN
    1 - \frac{\iota_k}{\cos^2\phi_k} + \ln\(\frac{\iota_k}{\cos^2\phi_k}\) > -c_0.
  \eequationNN
  Hence, by \eqref{eq.zeta}, there exist $c_2 \in \R{}_{>0}$ and $c_3 \in \R{}_{>0}$ such that, for all $k \in J_{p,K}$,
  \bequationNN
    c_2 \leq \frac{\iota_k}{\cos^2\phi_k} \leq c_3.
  \eequationNN
  Combining this with the fact (already proved) that $\cos\phi_k > c_1$ for all $k \in J_{p,K}$ and the fact that $\cos\phi_k \leq 1$ for all $k \in \N{}$, it follows, for all $k \in J_{p,K}$, that $c_1^2 c_2 \leq \iota_k \leq c_3$.  Therefore, since $\|H_ks_k\|_2/\|s_k\|_2 = \iota_k/\cos\phi_k$,
  \bequationNN
    c_1^2c_2 
    \leq \frac{\iota_k}{\cos\phi_k}
    = \frac{\|H_ks_k\|_2}{\|s_k\|_2}
    = \frac{\iota_k}{\cos\phi_k}
    \leq \frac{c_3}{c_1};
  \eequationNN
  i.e., the latter inequalities in \eqref{eq.good_bounds_H} hold for any $\sigma \in (0,c_1^2c_2]$ and $\mu \in [c^{-1}c_3,\infty)$.
\eproof

%*********
% Section
%*********
\section{Proofs for \ref{alg.bundle}}\label{app.bundle}

This appendix provides proofs for the results stated in \S\ref{subsec.bundle} related to the \ref{alg.bundle} algorithm.  The results and proofs are based on those found in \S7.4 in \cite{Rusz06}, but modified to account for a variable-metric quadratic term in the subproblem objective.

\bproof[Proof of Lemma~\ref{lem.MY_bound_f}]
  Let $\xtilde \in \Xcal_k$ be any point with $f(\xtilde) < f(x_k)$, the existence of which follows under the conditions of the lemma.  Restricting the minimization on the right-hand side of \eqref{eq.MY} to the line segment $[x_k,\xtilde]$ and using convexity of $f$ gives
  \bequationNN
    \baligned
    f_{H_k,\Xcal_k}(x_k)
      &\leq \min_{x \in [x_k,\xtilde]} \ f(x) + \thalf(x-x_k)^TH_k(x-x_k) \\
      &=    \min_{\Delta\in[0,1]} \ f((1-\Delta)x_k + \Delta\xtilde) + \Delta^2 \thalf (\xtilde - x_k)^TH_k(\xtilde - x_k) \\
      &\leq \min_{\Delta\in[0,1]} \ (1-\Delta)f(x_k) + \Delta f(\xtilde) + \Delta^2 \thalf (\xtilde - x_k)^TH_k(\xtilde - x_k) \\
      &=    f(x_k) + \min_{\Delta\in[0,1]} \ \Delta(f(\xtilde) - f(x_k)) + \Delta^2 \thalf (\xtilde - x_k)^TH_k(\xtilde - x_k).
    \ealigned
  \eequationNN
  Since $H_k \succ 0$, this last minimization over $\Delta \in [0,1]$ involves a strongly convex quadratic function of~$\Delta$.  Moreover, since $f(\xtilde) < f(x_k)$, the value of $\Delta \in [0,1]$ that minimizes the function is strictly positive, meaning that the optimal value of the problem is strictly negative.
\eproof

\bproof[Proof of Lemma~\ref{lem.MY_bound_q}]
  Let $(k,m) \in \N{} \times \N{}$ be given.  Since $H_k \succ 0$, it follows by \eqref{eq.q} that $l_{k,m}(x) \leq q_{k,m}(x)$ for all $x \in \R{n}$, giving the first inequality in \eqref{eq.lqf}.  Moreover, $l_{k,m}$ being a pointwise underestimator of $f$ throughout $\R{n}$ means that $l_{k,m}(x) \leq f(x)$ for all $x \in \R{n}$, which implies by \eqref{eq.q} that $q_{k,m}(x) \leq f(x) + \thalf (x-x_k)^TH_k(x-x_k)$ for all $x \in \R{n}$.  Letting $\xhat$ be the argument that solves the minimization problem in the right-hand-side of~\eqref{eq.MY} for $\xbar = x_k$, it follows along with the arguments above that
\bequationNN
  q_{k,m}(x_{k,m+1}) \leq q_{k,m}(\xhat) \leq f(\xhat) + \thalf (\xhat-x_k)^TH_k(\xhat-x_k) 
= f_{H_k,\Xcal_k}(x_k),
\eequationNN
which establishes the second inequality in~\eqref{eq.lqf}.
\eproof

\bproof[Proof of Lemma~\ref{lem.bundle.descent}]
  Let $k \in \N{}$ be given and suppose that $x_k$ is not a minimizer of~$f$.  Then, for any $m \in \N{}$, Lemmas~\ref{lem.MY_bound_f} and \ref{lem.MY_bound_q} imply that
  \bequation\label{eq.string}
    l_{k,m}(x_{k,m+1}) \leq q_{k,m}(x_{k,m+1}) \leq f_{H_k,\Xcal_k}(x_k) < f(x_k),
  \eequation
  meaning that the termination check in Step~\ref{step.termination_check} never tests true.  Hence, to derive a contradiction to the statement of the lemma, suppose that the algorithm generates an infinite sequence $\{x_{k,m+1}\}_{m=1}^\infty$ such that no element satisfies the condition in Step~\ref{step.decrease}.
  
  Toward deriving the aforementioned contradiction, let us first show that the generated function values $\{f(x_{k,m+1})\}_{m=1}^\infty$ converge to the minimizer of $f$ over~$\Xcal_k$, namely, $f_{\Xcal_k} := \min_{x\in\Xcal_k} f(x)$.  Notice that since $x_k$ is not a minimizer of $f$, it follows that $f_{\Xcal_k} < f(x_k)$.  For any $\varepsilon \in (0,\infty)$, let
  \bequationNN
    \Mcal_\varepsilon := \{m \in \N{} : f_{\Xcal_k} + \varepsilon < f(x_{k,m+1})\}.
  \eequationNN
  Suppose that there exists a pair $(m_1,m_2) \in \Mcal_\varepsilon \times \Mcal_\varepsilon$ with $m_1 < m_2$.  Then, since $\{l_{k,m}\}_{m=1}^\infty$ are pointwise underestimators of $f$, we can conclude that
  \bequation\label{eq.airplane}
    f_{k,m_1} + g_{k,m_1}^T(x_{k,m_2+1} - x_{k,m_1}) \leq l_{k,m_2}(x_{k,m_2+1}) \leq f_{\Xcal_k}.
  \eequation
  On the other hand, by virtue of $m_2$ being an element of $\Mcal_\varepsilon$, it follows that $\varepsilon < f(x_{k,m_2+1}) - f_{\Xcal_k}$, which combined with \eqref{eq.airplane} implies that
  \bequation\label{eq.eps}
    \varepsilon < f(x_{k,m_2+1}) - f_{k,m_1} - g_{k,m_1}^T(x_{k,m_2+1} - x_{k,m_1}).
  \eequation
  Since $\Xcal_k$ is compact, there exists $L_{\Xcal_k} \in (0,\infty)$ such that (recall~\eqref{eq.Lipschitz})
  \bequationNN
    |f(x) - f(\xbar)| \leq L_{\Xcal_k}\|x - \xbar\|_2\ \ \text{for all}\ \ (x,\xbar) \in \Xcal_k \times \Xcal_k.
  \eequationNN
  Since the subgradients of $f$ on $\Xcal_k$ are bounded, one can assume that $L_{\Xcal_k}$ is large enough such that $\|g_{k,m+1}\|_2 \leq L_{\Xcal_k}$ for all $m \in \N{}$.  Hence, from \eqref{eq.eps}, one finds
  \bequationNN
    \varepsilon < 2L_{\Xcal_k}\|x_{k,m_2+1} - x_{k,m_1}\|_2\ \ \text{for all}\ \ (m_1,m_2) \in \Mcal_\varepsilon \times \Mcal_\varepsilon\ \ \text{with}\ \ m_1 \neq m_2.
  \eequationNN
  Since the set $\Xcal_k$ is compact, there can only exist a finite number of points in $\Xcal_k$ having a distance at least~$\varepsilon/(2L_{\Xcal_k})$ from each other. Thus, the set $\Mcal_\varepsilon$ must be finite. In turn, this means that for any $\varepsilon \in (0,\infty)$ there can only be a finite number of points with objective function value in $[f_{\Xcal_k},f_{\Xcal_k} + \varepsilon]$. One may thus conclude that the sequence $\{f(x_{k,m+1})\}_{m=1}^\infty$ converges to $f_{\Xcal_k}$.
  
  Let us now use the established convergence of $\{f(x_{k,m+1})\}_{m=1}^\infty$ to $f_{\Xcal_k}$ to derive a contradiction to the supposition that $\{x_{k,m+1}\}_{m=1}^\infty$ is generated with no element satisfying the condition in Step~\ref{step.decrease}.  Since $\{x_{k,m+1}\}_{m=1}^\infty$ is contained in the compact set~$\Xcal_k$, there exists an infinite $\Mcal \subseteq \N{}$ such that
  \bequationNN
    \lim_{m\in\Mcal,m\to\infty} x_{k,m+1} = \xbar\ \ \text{for some}\ \ \xbar \in \Xcal_k.
  \eequationNN
  Since $\{f(x_{k,m+1})\}_{m=1}^\infty \to f_{\Xcal_k}$, it follows that $f(\xbar) = f_{\Xcal_k}$.  For any $m\in\Mcal$, let $\mbar$ be the smallest element in $\Mcal$ that is strictly larger than $m$.  It follows using the same argument that led to~\eqref{eq.airplane} that
\bequationNN
  f(x_{k,m+1}) + g_{k,m+1}^T(x_{k,\mbar+1} - x_{k,m+1}) \leq l_{l,\mbar}(x_{k,\mbar+1}) \leq f_{\Xcal_k}.
\eequationNN
  Taking limits over $m\in\Mcal$ as $m\to\infty$, using the uniform bound on the subgradient norms $\{\|g_{k,m+1}\|\}_{m=1}^\infty$ over $\Xcal_k$ as described in the previous paragraph, and recalling the facts that $\lim_{m\in\Mcal,m\to\infty} x_{k,m+1} = \xbar$ and $f(\xbar) = f_{\Xcal_k}$, one finds that
  \bequationNN
    f_{\Xcal_k} = f(\xbar) = \lim_{m\in\Mcal,m\to\infty} \(f(x_{k,m+1}) + g_{k,m+1}^T(x_{k,\mbar+1} - x_{k,m+1})\) \leq \lim_{m\in\Mcal,m\to\infty} l_{k,\mbar}(x_{k,\mbar+1}) \leq f_{\Xcal_k}.
  \eequationNN
  Since $\lim_{m\in\Mcal,m\to\infty} l_{k,m}(x_{k,m+1}) = \lim_{m\in\Mcal,m\to\infty} l_{k,\mbar}(x_{k,\mbar+1})$, this proves that
  \bequationNN
    \lim_{m\in\Mcal,m\to\infty} l_{k,m}(x_{k,m+1}) = f_{\Xcal_k},
  \eequationNN
  from which it follows that
  \bequationNN
    \lim_{m\in\Mcal,m\to\infty} \(\frac{f(x_k) - f(x_{k,m+1})}{f(x_k) - l_{k,m}(x_{k,m+1})}\) = \frac{f(x_k) - f_{\Xcal_k}}{f(x_k) - f_{\Xcal_k}} = 1.
  \eequationNN
  We have reached a contradiction since this limit indicates that the condition in Step~\ref{step.decrease} would be satisfied for some sufficiently large $m \in \Mcal$.
\eproof

%**************
% Bibliography
%**************
\bibliography{references}

\end{document}